\documentclass[12pt]{amsart}
\usepackage{amsmath,amscd,amssymb,amsfonts,graphics,mathrsfs}
\newcommand{\h}{\hbox}
\newcommand{\q}{\quad}
\newcommand{\nin}{\noindent}
\newcommand{\bs}{\par\bigskip}
\newcommand{\ms}{\par\medskip}
\newcommand{\sk}{\par\smallskip}
\newcommand{\bsn}{\par\bigskip\noindent}
\newcommand{\msn}{\par\medskip\noindent}
\newcommand{\skn}{\par\smallskip\noindent}
\newcommand{\ges}{\geqslant}
\newcommand{\les}{\leqslant}
\newcommand{\1}{\hskip1pt}
\newcommand{\mcap}{\hbox{$\bigcap$}}
\newcommand{\mcup}{\hbox{$\bigcup$}}

\newcommand{\msum}{\hbox{$\sum$}}
\newcommand{\mopl}{\hbox{$\bigoplus$}}
\newcommand{\mprod}{\hbox{$\prod$}}
\newcommand{\mwdg}{\hbox{$\bigwedge$}}
\newcommand{\Cc}{{\mathscr C}}
\newcommand{\D}{{\mathscr D}}

\newcommand{\F}{{\mathscr F}}

\newcommand{\Hc}{{\mathscr H}}
\newcommand{\I}{{\mathscr I}}
\newcommand{\J}{{\mathscr J}}
\newcommand{\K}{{\mathscr K}}
\newcommand{\Lc}{{\mathscr L}}
\newcommand{\Nc}{{\mathscr N}}
\newcommand{\M}{{\mathscr M}}

\newcommand{\OO}{{\mathscr O}}
\newcommand{\PP}{{\mathbb P}}
\newcommand{\Q}{{\mathbb Q}}
\newcommand{\C}{{\mathbb C}}
\newcommand{\N}{{\mathbb N}}
\newcommand{\R}{{\mathbb R}}
\newcommand{\RR}{{\mathbf R}}
\newcommand{\DD}{{\mathbf D}}
\newcommand{\Z}{{\mathbb Z}}

\newcommand{\Dsf}{{\sf D}}
\newcommand{\Msf}{{\sf M}}
\newcommand{\al}{\alpha}
\newcommand{\de}{\delta}
\newcommand{\De}{\Delta}
\newcommand{\ep}{\varepsilon}
\newcommand{\ga}{\gamma}
\newcommand{\Ga}{\Gamma}
\newcommand{\la}{\lambda}
\newcommand{\om}{\omega}
\newcommand{\Om}{\Omega}

\newcommand{\ft}{\widetilde{f}}

\newcommand{\Dt}{\widetilde{D}}

\newcommand{\St}{\widetilde{S}}

\newcommand{\Vt}{\widetilde{V}}

\newcommand{\Yt}{\widetilde{Y}}
\newcommand{\Zt}{\widetilde{Z}}

\newcommand{\fb}{{}\,\overline{\!f}{}}

\newcommand{\Ub}{\overline{U}}
\newcommand{\Xb}{{}\,\overline{\!X}{}}
\newcommand{\ub}{^{\ssb}}
\newcommand{\uD}{^{\D}}
\newcommand{\an}{^{\rm an}}

\newcommand{\IX}{\I_{\!X}}
\newcommand{\OX}{\OO_{\!X}}
\newcommand{\OZ}{\OO_{\!Z}}
\newcommand{\Diff}{{\rm Diff}}
\newcommand{\rh}{{\rm rh}}
\newcommand{\rhnc}{{\rm rh.nc}}
\newcommand{\DR}{{\rm DR}}
\newcommand{\DRt}{{\rm DR}^{\rm top}}
\newcommand{\Gr}{{\rm Gr}}
\newcommand{\dd}{\partial}
\newcommand{\ddd}{{\rm d}}

\newcommand{\ob}{{\mathbf 1}}
\newcommand{\sot}{{\otimes}}
\newcommand{\Hom}{{\mathscr H}\hskip-1.7pt om}
\newcommand{\Ext}{{\mathscr E}\hskip-1.4pt xt}
\newcommand{\lYr}{\langle Y\rangle}
\newcommand{\bl}{\bigl}
\newcommand{\br}{\bigr}
\newcommand{\pl}{\1{+}\1}
\newcommand{\mi}{\1{-}\1}
\newcommand{\eq}{\,{=}\,}
\newcommand{\stm}{\,{\setminus}\,}
\newcommand{\sst}{\,{\subset}\,}
\newcommand{\tos}{\,{\to}\,}
\newcommand{\caps}{\,{\cap}\,}

\newcommand{\ins}{\,{\in}\,}
\newcommand{\less}{\,{\leqslant}\,}
\newcommand{\gess}{\,{\geqslant}\,}

\newcommand{\ssb}{\raise.15ex\h{${\scriptscriptstyle\bullet}$}}
\newcommand{\ssc}{\,\raise.15ex\h{${\scriptstyle\circ}$}\,}
\newcommand{\onto}{\twoheadrightarrow}
\newcommand{\into}{\hookrightarrow}
\newcommand{\simto}{\,\,\rlap{\hskip1.5mm\raise1.4mm\hbox{$\sim$}}\hbox{$\longrightarrow$}\,\,}
\newcommand{\simot}{\,\,\rlap{\hskip1.9mm\raise1.4mm\hbox{$\sim$}}\hbox{$\longleftarrow$}\,\,}
\newcommand{\indlim}{\rlap{\raise-5.7pt\h{$\,\to$}}{\rm lim}}
\newcommand{\Mt}{\rlap{\hskip3.2mm\raise2.5pt\hbox{$\widetilde{\hbox{\,\,\,}}$}}{\mathscr M}}
\newcommand{\hs}{\hskip1.3mm}
\def\cond[#1]{\par\noindent\rlap{#1}\hskip.7cm\hangindent=.7cm\hangafter=1}
\begin{document}
\h{}\bs
\centerline{\large Notes on Regular Holonomic $\D$-modules for Algebraic Geometers}
\bs
\centerline{Morihiko Saito}
\ms\bsn
\vbox{\narrower\smaller\nin{\bf Abstract.} We explain a formalism of regular holonomic $\D$-modules for algebraic geometers using the distinguished triangles associated with algebraic local cohomology together with meromorphic Deligne extensions of local systems as well as the dual functor.}
\ms\bs
\centerline{\bf Introduction}
\bsn
It has been noticed that regular holonomic $\D$-modules can be defined by induction on the dimension of support using distinguished triangles associated to {\it algebraic local cohomology\1} together with {\it meromorphic Deligne extensions\1} of local systems \cite[Prop.\,II.5.4]{De}, see \cite[B.5]{ypg}, \cite[Rem.\,1.2\,(ii)]{lcd}, etc. This approach may be easier to understand for algebraic geometers, although there are certain technical difficulties; for instance, the {\it stability by quotients\1} in the category of coherent $\D$-modules is not necessarily easy to show, and the {\it stability by proper direct images\1} must be proved at the same time by induction, where the {\it dual functor\1} $\DD$ is also needed.
\sk
In these notes we explain a formalism of regular holonomic $\D$-modules along this line, and show the following.
\msn
{\bf Theorem~1.} {\it For any complex algebraic variety or analytic space $X$, there is an abelian category $\Msf_{\rh}(\D_X)$ of regular holonomic $\D_X$-modules which are stable by the dual functor $\DD$ and by subquotients and extensions in the category of coherent $\D_X$-modules. Moreover the categories $\Msf_{\rh}(\D_X)$ are stable by the functors $\Hc^j_{[Y]}$, $\Hc^j_{[X|Y]}$ for any closed subvarieties $Y\sst X$ and by the direct images $\Hc^j\!f^{\D}_*$ and the pullbacks $\Hc^j\!f_{\D}^!$ for any morphisms $f$ except for the direct images in the analytic case where $f$ is assumed to be proper.}
\ms
For proper direct images, we extend the {\it calculation of nearby cycles\1} in \cite{St}, \cite{gau} by using the {\it Hartogs extension theorem\1} (see Rem.\,1.3f below) to simplify some arguments (where the relation to the calculation of Gauss-Manin connections in \cite{De} does not seem trivial). One advantage of this method is that we can quite easily prove {\it constructibility of de Rham complexes\1} of regular holonomic $\D$-modules by reducing to the normal crossing case. (For other ways, see \cite{Be}, \cite{Bo}, \cite{HT}, \cite{KK}, \cite{Me2}, etc.)
\sk
In Section 1 we first recall some basics of holonomic $\D$-modules in {\bf 1.1} and algebraic local cohomology functor in {\bf 1.2} assuming some elementary $\D$-module theory which is partly explained in Appendix. We then review regular holonomic $\D$-modules of normal crossing type including Deligne extensions in {\bf 1.3} together with regular meromorphic connections in {\bf 1.4} and relative de Rham complexes over a disk in {\bf 1.5}. In Section 2 we give a definition of regular holonomic $\D$-modules using algebraic local cohomology and meromorphic Deligne extensions in {\bf 2.1}, and prove some of their properties in {\bf 2.2}. In Appendix, we review some elementary $\D$-module theory.
\sk
This work was partially supported by JSPS Kakenhi 15K04816.
\msn
{\bf Convention~1.} In these notes, a {\it complex variety\1} means either a {\it complex algebraic variety\1} in the sense of Serre \cite[\S 34]{Se} or a {\it reduced complex analytic space.} In the algebraic case, it is essentially a separated reduced scheme of finite type over $\C$ except that we consider only {\it closed points.} Here {\it irreducibility\1} of varieties is not supposed, and $\OX$-modules in the algebraic case are {\it not\1} necessarily quasi-coherent (in the usual sense).
\sk
{\bf 2.} An $\OX$-module $\M$ is called {\it quasi-coherent\1} in these notes if there is locally an increasing exhaustive filtration $G$ of $\M$ indexed by $\N$ such that the $G_k\M$ are {\it coherent\1} $\OX$-submodules for any $k\in\Z$. This is suitable in the case of a coherent $\D_X$ (or $\OX(*Y)$)-module $\M$ (with $Y\sst X$ a divisor), since $\M$ is coherent over $\D_X$ (or $\OX(*Y)$) if and only if it is locally finitely generated over $\D_X$ (or $\OX(*Y)$), and is quasi-coherent over $\OX$ in the above sense.
\bs\bs
\vbox{\centerline{\bf 1. Preliminaries on holonomic $\D$-modules}
\bsn
In this section we first recall some basics of holonomic $\D$-modules in {\bf 1.1} and algebraic local cohomology functor in {\bf 1.2} assuming some elementary $\D$-module theory which is partly explained in Appendix. We then review regular holonomic $\D$-modules of normal crossing type including Deligne extensions in {\bf 1.3} together with regular meromorphic connections in {\bf 1.4} and relative de Rham complexes over a disk in {\bf 1.5}.}
\msn
{\bf 1.1.~Holonomic $\D$-modules.} Let $\D_X$ be the sheaf of differential operators on a smooth complex variety $X$. It has the filtration $F$ by the order of differential operators, and its associated graded ring $\Gr^F_{\ssb}\D_X$ is locally isomorphic to the polynomial ring $\OX[\xi_1,\dots,\xi_{d_X}]$, where $\xi_i:=\Gr^F_1\dd_{z_i}$ with $z_1,\dots,z_{d_X}$ local coordinates. In the case of algebraic varieties, local coordinates give only an \'etale morphism to $\C^{d_X}$; however, the associated vector fields $\dd_{z_i}$ are defined by the condition $\langle\dd_{z_i},\ddd z_j\rangle\eq\de_{i,j}$ using $\Theta_X\eq\Hom_{\OX}(\Om_X^1,\OX)$ as is well-known. (In these notes, $\Theta_X$ is the sheaf of vector fields.)
\sk
In the algebraic case, ${\rm Spec}_X\Gr^F_{\ssb}\D_X$ is then naturally identified with the cotangent bundle $T^*\!X$. In the analytic case, ${\rm Spec}_X\Gr^F_{\ssb}\D_X$ can be replaced with ${\rm Specan}_X\Gr^F_{\ssb}\D_X$. This means that $\Gr^F_{\ssb}\D_X$ is identified with the direct image to $X$ of the sheaf of holomorphic functions on $T^*\!X$ whose restrictions to the fibers are polynomials.
\sk
We say that a coherent left $\D_X$-module $\M$ is {\it holonomic\1} if $M$ has locally a coherent filtration $F$ such that ${\rm Supp}\,\Gr^F_{\ssb}\M\sst T^*\!X$ (which is called the {\it characteristic variety\1} of $\M$, and is denoted by ${\rm Char}\,\M$) has dimension $d_X\,{:=}\,\dim X$. Here a {\it coherent filtration\1} means that $(\M,F)$ is a filtered $\D_X$-module, $F_p\1\M\eq0$ for $p\,{\ll}\,0$ (locally on $X$), and $\Gr^F_{\ssb}\M$ is a coherent $\Gr^F_{\ssb}\D_X$-module.
\sk
We denote the category of holonomic left $\D_X$-modules by $\Msf_{\rm hol}(\D_X)$. This is an abelian subcategory of the category of left $\D_X$-modules $\Msf(\D_X)$, and is closed by subquotients and extensions in $\Msf(\D_X)$. It is stable by the dual functor $\DD$ (see Rem.\,1.1 below) which can be defined locally for a coherent left $\D_X$-module $\M$ by
$$\DD(\M)=\RR\Hom_{\D_X}(\M,\D_X\sot_{\OX}\om_X^{\vee})[d_X],
\leqno(1.1.1)$$
taking locally a free resolution of $\M$, see also {\bf A.3} below. Here $\D_X\sot_{\OX}\om_X^{\vee}$ means that the right $\D_X$-module structure of $\D_X$ is changed to the left one, see {\bf A.1} below. (More precisely, it is better to note it as $\Hom_{\OX}(\om_X,\D_X)$, although it does not seem very clear to the reader which $\D_X$-module structure of $\D_X$ is used for $\Hom_{\OX}$ without indicating it explicitly in this case.)
\sk
We will denote by $\Dsf^b(\D_X)$ the derived category of bounded complexes of left $\D_X$-module, and by $\Dsf^b_{\rm coh}(\D_X)$, $\Dsf^b_{\rm hol}(\D_X)$ its full subcategories consisting of complexes having respectively coherent and holonomic $\D_X$-modules as cohomology sheaves.
(In the algebraic case, we do {\it not\1} assume that $\D_X$-modules are quasi-coherent over $\OX$.)
\msn
{\bf Remark\hs1.1.} It is well-known that for $\M\ins \Msf_{\rm hol}(\D_X)$, we have
$$\DD(\M)\ins \Msf_{\rm hol}(\D_X)\,\,\,\,\h{(in particular}\,\,\,\Hc^j\DD(\M)\eq0\,\,\,\h{if}\,\,\,j\,{\ne}\,0).
\leqno(1.1.2)$$
This follows easily from the {\it involutivity\1} of characteristic variety implying that its dimension is at least $d_X$, see \cite{SKK}, \cite{Mal}, etc. Taking locally a filtered free resolution of $(\M,F)$, the dual functor $\DD$ is compatible with the passage to the associated graded quotients of a coherent filtration $\Gr^F_{\ssb}\M$. This gives the vanishing of $\Hc^j\DD(\M)$ ($j\,{\ne}\, 0$) at {\it generic\1} points of characteristic variety, hence everywhere on $T^*\!X$ by the {\it involutivity,} see also \cite{Bj}, \cite{KaB}, \cite{Ka4}, \cite{Sab}, etc. Indeed, setting $(\M'{}\ub,F):=\DD(\M,F)$ locally on $X$, we have the spectral sequence
$$E_1^{p,q}=\Hc^{p+q}\1\Gr^F_{-p}\1\M'{}\ub\Longrightarrow\Hc^{p+q}\M'{}\ub,
\leqno(1.1.3)$$
which is compatible with the action of any vector field $\xi$, and $\Gr^F_{\ssb}\Hc^j\!\M'{}\ub$ is a {\it subquotient of $\Hc^j\Gr^F_{\ssb}\M'{}\ub$ as a graded $\Gr^F_{\ssb}\D_X$-module.} Here we can also consider the complex $\M'{}\ub[t,t^{-1}]$ with the filtration $F_{\ssb}$ defined by $\mopl_{k\in\Z}\,F_{\ssb+k}\1\M'{}\ub\1t^k$ so that it is a complex of filtered modules over the graded ring $\mopl_{k\in\Z}\,F_k\D_X\1t^k$, and its $p$\1th graded piece is $\mopl_{k\in\Z}\,\Gr^F_{p+k}\1\M'{}\ub\1t^k$. (It does not seem clear if there is always a coherent filtration $F$ on any holonomic $\D_X$-module $\M$ locally on $X$ so that the $\Gr^F_{\ssb}\M_x$ are {\it Cohen-Macaulay\1} $\Gr^F_{\ssb}\D_{X,x}$-modules.)
\msn
{\bf 1.2.~Algebraic local cohomology.} Let $X$ be a smooth complex variety. Let $Y\sst X$ be a closed subvariety with $\I_Y$ its ideal sheaf. For an $\OX$-module $\M$, the {\it algebraic local cohomology\1} $\Hc^j_{[Y]}\M$ is defined by
$$\Hc^j_{[Y]}\M:=\rlap{\raise-10pt\h{$\,\,\,\scriptstyle k$}}\indlim\,\Ext^j_{\OX}(\OX/\!\1\I_Y^k,\M)\q\q\q(j\ins\Z),
\leqno(1.2.1)$$
see \cite{Gr2}. Here $\I_Y$ can be replaced by any coherent ideal sheaf $\I$ with ${\rm Supp}\,\OX/\!\1\I=Y$. We can extend this to $\D_X$-modules $\M$ using injective resolutions. Note that
$$\rlap{\raise-10pt\h{$\,\,\,\scriptstyle k$}}\indlim\,\Hom_{\OX}(\OX/\!\1\I_Y^k,\M)=\mcup_k\,\M^{\I_Y^k},
\leqno(1.2.2)$$
where the right-hand side is the union of the subsheaf of $\M$ annihilated by $\I_Y^k$. Indeed, we have $g(\dd_{z_i}m)\eq\dd_{z_i}(gm)\mi(\dd_{z_i}g)m\eq0$ for $g\ins\I_Y^{k+1}$, $m\ins\M^{\I_Y^k}$, since $\dd_{z_i}g\ins\I_Y^k$. Here the $z_i$ are local coordinates, see also \cite{Ka2}. (Recall that injective $\D_X$-modules are injective $\OX$-modules, since $\D_X$ is flat over $\OX$.)
\sk
Similarly we can define the functors
$$\RR\Ga_{[Y]},\q\RR\Ga_{[X|Y]},\q\Hc^j_{[X|Y]},
\leqno(1.2.3)$$
replacing $\Ext^j_{\OX}(\OX/\!\1\I_Y^k,\M)$ in the right-hand side of (1.2.1) respectively with
$$\Hom_{\OX}(\OX/\!\1\I_Y^k,\J\ub),\q\Hom_{\OX}(\I_Y^k,\J\ub),\q\Ext^j_{\OX}(\I_Y^k,\M),$$
where $\J\ub$ is an injective resolution of $\M$ as a $\D_X$-module.
(It may intuitively easier to understand the notation if we write $\scriptstyle[X{\setminus}Y]$ rather than $\scriptstyle[X|Y]$. However, this may cause a problem when $X{\setminus}Y$ is a {\it closed\1} subvariety of some smooth complex variety.)
\sk
The first two functors of (1.2.3) is naturally extended to functors from $\Dsf^b(\D_X)$ to $\Dsf^b(\D_X)$ using the truncation $\tau_{\les k}$ for $k\gg 0$ (since the $\OO_{X,x}$ are regular local rings), see also (A.2.8) below. (Recall that $\tau_{\les k}\M\ub$ is defined by $(\tau_{\les k}\M\ub)^j\eq\M^j$ if $j\,{<}\,k$, ${\rm Ker}(\ddd\,{:}\,\M^k\tos\M^{k+1})$ if $j\eq k$, and $0$ otherwise.)
\sk
In the algebraic case, $\RR\Ga_{[Y]},\Hc^j_{[Y]}$ are usually denoted as $\RR\Ga_Y,\Hc^j_Y$, and we have the equality $\RR\Ga_{[X|Y]}\eq\RR(j_{X\setminus Y})_*j_{X\setminus Y}^{-1}$ for quasi-coherent sheaves with $j_{X\setminus Y}:X\stm Y\into X$ the inclusion.
\sk
In the $Y$ divisor case, we have
$$\RR\Ga_{[X|Y]}=\Ga_{[X|Y]}:=\Hc^0_{[X|Y]},
\leqno(1.2.4)$$
and it coincides with the usual localization along $Y$. Note that we have in general
$$\Hc^j_{[Y]}=\Hc^j\RR\Ga_{[Y]},\q\Hc^j_{[X|Y]}=\Hc^j\RR\Ga_{[X|Y]}.
\leqno(1.2.5)$$
\sk
By definition there is a distinguished triangle
$$\RR\Ga_{[Y]}\to{\rm id}\to\RR\Ga_{[X|Y]}\buildrel{+1}\over\to.
\leqno(1.2.6)$$
\sk
For a closed subvariety $Y'\sst Y$, the natural inclusion $\I_Y\sst\I_{Y'}$ induces the canonical morphisms
$$\RR\Ga_{[Y']}\to\RR\Ga_{[Y']},\q\RR\Ga_{[X|Y']}\to\RR\Ga_{[X|Y]}.
\leqno(1.2.7)$$
\sk
We have the following.
\msn
{\bf Lemma~1.2a.} {\it Let $Y'\subset Y$ be a closed subvariety. Then
$$\RR\Ga_{[Y']}\RR\Ga_{[Y]}=\RR\Ga_{[Y']},
\leqno(1.2.8)$$
hence there is the distinguished triangle}
$$\RR\Ga_{[Y']}\to\RR\Ga_{[Y]}\to\RR\Ga_{[X|Y']}\RR\Ga_{[Y]}\buildrel{+1}\over\to.
\leqno(1.2.9)$$
\msn
{\it Proof.} This follows from the assertion that $\Hom_{\OX}(\OX/\!\1\I,\J)$ is an injective $(\OX/\!\1\I)$-module, if $\I\subset\OX$ is a coherent ideal and $\J$ is an injective $\OX$-module. For the proof of the last assertion, we have the canonical isomorphism
$${\rm Hom}_{\OX/\!\I}\bl(\F,\Hom_{\OX}(\OX/\!\1\I,\J)\br)={\rm Hom}_{\OX}(\F,\J),$$
for $(\OX/\!\1\I)$-modules $\F$. This finishes the proof of Lem.\,1.2a.
\ms
From \cite[Thm.\,A1.3]{Eis} (applied to the stalks), we can deduce the following.
\msn
{\bf Lemma~1.2b.} {\it Assume that $Y=\mcap_{i=1}^r\1Z_i$ with $Z_i\sst X$ divisors. For $I\sst\{1,\dots,r\}$, set $Z_I:=\mcup_{i\in I}\,Z_i$. For a $\D_X$-module $\M$, there are \v Cech-type complexes
$$\RR\Ga_{[(Z_1,\dots,Z_r)]}\M,\q\RR\Ga_{[X|(Z_1,\dots,Z_r)]}\M,$$
such that their $j$\1th components are respectively
$$\mopl_{|I|=j}\,\Ga_{[Z_I]}\M,\q\mopl_{|I|=j+1}\,\Ga_{[Z_I]}\M\q\h{for}\q j\ges 0,
\leqno(1.2.10)$$
and $0$ otherwise $($where $\Ga_{[Z_{\emptyset}]}\M:=\M)$, and there are canonical isomorphisms
$$\aligned\RR\Ga_{[(Z_1,\dots,Z_r)]}\M&{}\eq\RR\Ga_{[Y]}\M,\\ \RR\Ga_{[X|(Z_1,\dots,Z_r)]}\M&{}\eq\RR\Ga_{[X|Y]}\M.\endaligned
\leqno(1.2.11)$$
These are naturally extended to the case $\M\ub\ins \Dsf^b(\D_X)$.}
\ms
Using this, we can show the following.
\msn
{\bf Proposition\hs1.2.} {\it Let $\pi:X'\to X$ be a proper morphism of smooth complex varieties. Set $Y':=\pi^{-1}(Y)$ with $Y=\mcap_{i=1}^r\,Z_i$ and $Z_i\sst X$ divisors. Then we have the isomorphisms}
$$\aligned\RR\Ga_{[Y]}\ssc\pi_*\uD&=\pi_*\uD\ssc\RR\Ga_{[Y']},\\ \RR\Ga_{[X|Y]}\ssc\pi_*\uD&=\pi_*\uD\ssc\RR\Ga_{[X'|Y']}.\endaligned
\leqno(1.2.12)$$
\msn
{\it Proof.} The direct image can be defined by using the factorization
$$\pi:X'\into X'{\times}X\onto X,$$
where the first morphism is the graph embedding of $\pi$. The second one is the second projection, and the relative de Rham complex is used for the definition of direct image. Using the \v Cech-type complexes as in Lem.\,1.2b, we then see that the assertion is reduced to the second isomorphism of (1.2.12) in the $Y$ divisor case. It is then easy to verify it in the second projection case, since the sheaf-theoretic direct image by a proper morphism commutes with inductive limit. We may thus assume $\pi$ is a closed immersion $i:X'\into X$ and $Y\sst X$ is a divisor.
\sk
For a $\D_{X'}$-module $\M'$, we have the canonical morphism
$$\iota:i_*\uD\M'\to i_*\uD\M'(*Y'),$$
where $(*Y')$ denotes the usual localization along $Y'\eq X'\caps Y$. It is then enough to show that the action on $i_*\uD\M'(*Y')$ of a local defining function $g$ of $Y\sst X$ is bijective using the (exact) functoriality of localization. (Note that the kernel and cokernel of $\iota$ are quasi-coherent $\OX$-modules supported on $Y'\sst Y$, and are annihilated by the localization along $Y$.)
\sk
By definition $i_*\uD\M'(*Y')$ is locally isomorphic to $\M'(*Y')[\dd_{z_1},\dots,\dd_{z_r}]$, where $z_1,\dots,z_{d_X}$ are local coordinates of $X$ such that $X'\sst X$ is locally defined by $z_i\eq 0$ ($i\in[1,r])$. We have the expansion of the defining function
$$g=\msum_{\nu\in\Z^r}\,g_{\nu}\1z^{\nu}\q\bl(g_{\nu}\ins\C\{z_{r+1},\dots,z_{d_X}\},\,\,z^{\nu}\eq\mprod_{i=1}^r\,z_i^{\nu_i}\br).$$
Let $F$ be the increasing filtration of $\M'(*Y')[\dd_{z_1},\dots,\dd_{z_r}]$ by the order of $\dd_{z_i}$ ($i\in[1,r]$). Then $g_0$ is a defining function of $Y'\sst X'$ and the action of $g$ on $\Gr^F_{\ssb}\M'(*Y')[\dd_{z_1},\dots,\dd_{z_r}]$ coincides with that of $g_0$, since the action of $z_i$ ($i\in[1,r]$) on it vanishes. So the bijectivity of the action of $g$ follows. This finishes the proof of Prop.\,1.2.
\msn
{\bf Remark\hs1.2a.} For a coherent $\D_X$-module with ${\rm Supp}\,\M\sst Y$, we have $\RR\Ga_{[X|Y]}\M\eq 0$, that is, $\Hc^j_{[X|Y]}\M\eq 0$ ($\forall\,j\ins\Z$). This follows from Lem.\,1.2b, since the assertion is {\it local.}
This implies that, if there is a morphism of coherent $\D_X$-modules $\M'\tos\M$ such that the supports of kernel and cokernel are contained in $Y$, it induces the isomorphisms
$$\aligned\RR\Ga_{[X|Y]}\M'&\simto\RR\Ga_{[X|Y]}\M,\\ \Hc^j_{[X|Y]}\M'&\simto\Hc^j_{[X|Y]}\M\q(j\ins\Z).\endaligned$$
\msn
{\bf Remark\hs1.2b.} If there is $\M\ub\ins \Dsf^b_{\rm hol}(\D_X)$ with $\Hc^j\M\ub|_{X\setminus Y}\eq0$ ($j\ne0$) and $Y\sst X$ is a {\it divisor,} then $\Hc^j\RR\Ga_{[X|Y]}\M\ub\eq0$ ($j\ne 0$). Indeed, $\RR\Ga_{[X|Y]}\eq\Ga_{[X|Y]}$ is an exact functor (see (1.2.4)), and commutes with the passage to cohomology sheaves $\Hc^j$.
\msn
{\bf Remark\hs1.2c.} Assume that $Y=\mcap_{i=1}^r\,Z_i$ with $Z_i\sst X$ divisors. For a divisor $D\sst X$, set $Y':=Y\cup D$. Then we have the canonical isomorphism
$$\Hc^j_{[X|Y']}=\Hc^0_{[X|D]}\Hc^j_{[X|Y]}\q\q(j\ins\Z).$$
This follows from Lem.\,1.2b, since $\Hc^0_{[X|Z_I\cup D]}\eq\Hc^0_{[X|D]}\Hc^0_{[X|Z_I]}$ and $Y\cup D\eq\mcap_{i=1}^r(Z_i\cup D)$.
\msn
{\bf Remark\hs1.2d} ({\it Dual local cohomology functor\1}). For a holonomic left $\D_X$-module $\M$, set
$${}^{\vee}\!\!\Hc^j_{[Y]}\M:=\DD\Hc^{-j}_{[Y]}\DD\M,\q{}^{\vee}\!\!\Hc^j_{[X|Y]}\M:=\DD\Hc^{-j}_{[X|Y]}\DD\M.$$
There are exact sequences
$$0\to\Hc^0_{[Y]}\M\to\M\to\Hc^0_{[X|Y]}\M\to\Hc^1_{[Y]}\M\to0,$$
$$0\to{}^{\vee}\!\!\Hc^{-1}_{[Y]}\M\to{}^{\vee}\!\!\Hc^0_{[X|Y]}\M\to\M\to{}^{\vee}\!\!\Hc^0_{[Y]}\M\to0,$$
together with the commutative diagram
$$\begin{array}{cccccc}{}^{\vee}\!\!\Hc^0_{[X|Y]}\M&&\!\!\!\!\!\!\buildrel{\iota}\over\longrightarrow&&\!\!\Hc^0_{[X|Y]}\M\\ &\!\!\searrow&&\!\!\!\!\!\!\nearrow\\&&\!\!\!\!\!\!\M\\ &\!\!\nearrow&&\!\!\!\!\!\!\searrow\\\,\,\,\,\Hc^0_{[Y]}\M&&&&\!\!\!\!{}^{\vee}\!\!\Hc^0_{[Y]}\M
\end{array}$$
Note that $\Hc^0_{[Y]}\M$ (resp. ${}^{\vee}\!\!\Hc^0_{[Y]}\M$) is the largest sub (resp. quotient) object of $\M$ supported in $Y$. Hence ${\rm Im}({}^{\vee}\!\!\Hc^0_{[X|Y]}\M\tos\M)$ (resp. ${\rm Im}(\M\tos\Hc^0_{[X|Y]}\M)$) is the smallest sub (resp. quotient) object of $\M$ whose restriction to $X\stm Y$ coincides with that of $\M$.
\sk
In the $\M$ regular holonomic case (see {\bf 2.1} below), the image of $\iota$ is called the {\it minimal extension\1} of $\M|_{X\setminus Y}$. It has no non-trivial sub nor quotient object supported in $Y$. This can be shown by using the above commutative diagram for the image of $\iota$ (instead of $\M$).
\msn
{\bf 1.3.~Regular holonomic $\D$-modules of normal crossing type.} Let $X$ be a complex manifold with $z\eq(z_1,\dots,z_{d_X})$ a local coordinate system around $x\ins X$. For a holonomic left $\D_{X,x}$-module $\M_x$ and for $\nu=(\nu_1,\dots,\nu_{d_X})\ins\C^{d_X}$, set
$$\M_{x,z}^{(\nu)}:=\bl\{m\ins\M_x\mid(z_i\dd_{z_i}{-}\1\nu_i)^k\1m=0\,\,\,(k\gg0,\,\,\forall\,i\in[1,d_X])\br\}.$$
We say that $\M_x$ is {\it regular holonomic of normal crossing type,} if there is a local coordinate system $z=(z_1,\dots,z_{d_X})$ around $x$ such that $\M_x$ is generated over $\D_{X,x}$ by $\M_{x,z}^{(\nu)}$ for $\nu\ins\Xi^{d_X}$, where
$$\Xi:=\Xi'\cup\,\{-1\}\q\h{with}\q\Xi':=\bl\{\theta\ins\C\mid {\rm Re}\,\theta\in(-1,0]\br\}.$$
The category of such $\D_{X,x}$-modules with respect to a local coordinate system $z$ will be denoted by $M^{(z)}_{\rhnc}(\D_{X,x})$. (Note that $\M_{x,z}^{(\nu)}=0$ except for a finite number of $\nu\ins\Xi^{d_X}$, since $\M$ is holonomic.) We say that a left $\D_X$-module $\M$ is {\it regular holonomic of normal crossing type\1} if so is $\M_x$ for any $x\ins X$.
\msn
{\bf Remark\hs1.3a.} For $\nu=(\nu_1,\dots,\nu_{d_X})\ins\C^{d_X}$, we have the bijections
$$z_i:\M_{x,z}^{(\nu)}\simto\M_{x,z}^{(\nu+\ob_i)},\q\dd_{z_i}:\M_{x,z}^{(\nu+\ob_i)}\simto\M_{x,z}^{(\nu)}\q\h{if}\q\nu_i\ne -1.
\leqno(1.3.1)$$
(Here the $j$\1th component of $\ob_i\ins\Z^{d_X}$ is $1$ if $j\eq i$, and 0 otherwise.)
So it is enough to consider the $\M_{x,z}^{(\nu)}$ with $\nu\ins\Xi^{d_X}$ as in the above definition.
\msn
{\bf Remark\hs1.3b} ({\it Combinatorial description\1}). For $\M_x\ins M^{(z)}_{\rhnc}(\D_{X,x})$, we have the natural inclusions
$$\mopl_{\nu\in\C^{d_X}}\,\M_{x,z}^{(\nu)}\,\into\,\M_x\into\mprod_{\nu\in\C^{d_X}}\,\M_{x,z}^{(\nu)}.
\leqno(1.3.2)$$
Here the first term is dense in the middle term with respect to the canonical DFS (dual Fr\'echet-Schwartz) topology, and its completion by this topology (that is, the closure in $\M_x$) can be denoted by $\widehat{\mopl}_{\nu}\,\M_{x,z}^{(\nu)}$.
This coincides with the submodule of $\mprod_{\nu}\,\M_{x,z}^{(\nu)}$ generated by $\mopl_{\nu}\,\M_{x,z}^{(\nu)}$ over $\OO_{X,x}$. Combining this with Rem.\,1.3a, we see that $\M_x\ins M^{(z)}_{\rhnc}(\D_{X,x})$ is uniquely determined by the finite dimensional vector spaces $\M_{x,z}^{(\nu)}$ for $\nu\ins\Xi^{d_X}$ together with the morphisms $z_i\dd_{z_i}{-}\1\nu_i$, $z_i$, $\dd_{z_i}$ among them.
\msn
{\bf Remark\hs1.3c} ({\it Localization\1}). For $\M_x\ins M^{(z)}_{\rhnc}(\D_{X,x})$, we have $\Hc^0_{[X|Y_i]}\M_x\ins M^{(z)}_{\rhnc}(\D_{X,x})$ with $Y_i\eq\{z_i\eq 0\}$. Indeed, we can modify the $\M_{x,z}^{(\nu)}$ for $\nu\ins\Xi^{d_X}$ with $\{i\mid\nu_i\eq{-1}\}\ne\emptyset$ so that the morphisms $z_i$ in (1.3.1) give isomorphisms when $\nu_i\eq {-1}$.
\msn
{\bf Proposition\hs1.3a.} {\it The category $M^{(z)}_{\rhnc}(\D_{X,x})$ is stable by subquotients and extensions in the abelian category of $\D_{X,x}$-modules $\Msf(\D_{X,x})$. It is also stable by the dual functor $\DD$.}
\msn
{\it Proof.} Assume there is a short exact sequence of left $\D_{X,x}$-modules
$$0\to\M'_x\to\M_x\to\M''_x\to0.$$
Assume first $\M_x\ins M^{(z)}_{\rhnc}(\D_{X,x})$. By definition we have the exact sequences
$$0\to\M_{x,z}^{\prime\1(\nu)}\to\M_{x,z}^{(\nu)}\to\M_{x,z}^{\prime\prime\1(\nu)}\q\q(\nu\ins\C^{d_X}).$$
These imply the commutative diagram
$$\begin{array}{ccccccccc}0&\to&\widehat{\mopl}_{\nu}\M_{x,z}^{\prime\1(\nu)}&\to&\widehat{\mopl}_{\nu}\M_{x,z}^{(\nu)}&\to&\widehat{\mopl}_{\nu}\M_{x,z}^{\prime\prime\1(\nu)}\\ &&\downarrow&&\downarrow&&\downarrow&\raise5mm\h{}\raise-3mm\h{}\\ 0&\to&\M'_x&\to&\M_x&\to&\M''_x&\!\!\to&0\end{array}
\leqno(1.3.3)$$
where the middle vertical morphism is bijective. This implies the surjectivity of the right vertical morphism so that $\M''_x\ins M^{(z)}_{\rhnc}(\D_{X,x})$. The right vertical morphism is then bijective, and so is the left vertical morphism by the diagram, thus $\M'_x\ins M^{(z)}_{\rhnc}(\D_{X,x})$.
\sk
Assume now $\M'_x,\M''_x\ins M^{(z)}_{\rhnc}(\D_{X,x})$. Set
$$\M_{x,z_1}^{(\nu_1)}:=\bl\{m\ins\M_x\mid(z_1\dd_{z_1}{-}\1\nu_1)^k\1m=0\,\,\,(k\gg0)\br\},
\leqno(1.3.4)$$
and similarly for $\M_{x,z_1}^{\prime(\nu_1)},\M_{x,z_1}^{\prime\prime(\nu_1)}$. We have
$$\M_{x,z_1}^{\prime(\nu_1)},\M_{x,z_1}^{\prime\prime(\nu_1)}\ins M^{(z')}_{\rhnc}(\D_{X'{},x}),
\leqno(1.3.5)$$
where $X':=\{z_1\eq 0\}\sst X$ and $z':=(z_2,\dots,z_{d_X})$. By induction on $d_X$, it is then enough to show the short exact sequences
$$0\to\M_{x,z_1}^{\prime(\nu_1)}\to\M_{x,z_1}^{(\nu_1)}\to\M_{x,z_1}^{\prime\prime(\nu_1)}\to 0\q\q(\forall\,\nu_1\ins\C),
\leqno(1.3.6)$$
using a diagram similar to (1.3.3). It is easy to show the exactness except for the surjectivity of the last morphism. Take $m''\ins\M_{x,z_1}^{\prime\prime(\nu_1)}$. Let $m\ins\M_x$ whose image in $\M''_x$ is $m''$. There is $k\gg 0$ such that
$$(z_1\dd_{z_1}{-}\1\nu_1)^k\1m\ins\M'_x.$$
Since $(z_1\dd_{z_1}{-}\1\nu_1)^k$ is bijective on $\M_{x,z_1}^{\prime(\nu'_1)}$ if $\nu'_1\ne\nu_1$, we can modify $m$ by adding an element of $\M'_x$ so that
$$(z_1\dd_{z_1}{-}\1\nu_1)^k\1m\in\M_{x,z_1}^{\prime(\nu_1)}.$$
(Here we have to verify a convergence.) We then get that $(z_1\dd_{z_1}{-}\1\nu_1)^k\1m=0$ replacing $k$. So the first assertion follows. 
\sk
To show the last assertion, we may assume that $\M_x$ is simple, that is, $\M_{x,z}^{(\nu)}\cong\C$ for some $\nu\in\Xi^n$, and $\M_{x,z}^{(\nu')}\eq0$ for any $\nu\1'\in\Xi^{d_X}\stm\{\nu\}$. (This can be verified by induction on $d_X$ using Rem.\,1.3d just below.) Forgetting the $i\in[1,d_X]$ with $\nu_i\eq{-1}$ (and changing $d_X$), we may assume that $\nu_i\ne{-1}$ for any $i$. (This corresponds to taking the pullback by the closed embedding, see (A.4.3).) Then the assertion is easily verified. This finishes the proof of Prop.\,1.3a.
\msn
{\bf Remark\hs1.3d} ({\it Inductive description\1}). There is an inductive description of the category $\Msf_{\rhnc}^{(z)}(\D_{X,x})$ as follows. For $\M_x\ins \Msf_{\rhnc}^{(z)}(\D_{X,x})$, define $\M_{x,z_1}^{(\nu_1)}\ins M^{(z')}_{\rhnc}(\D_{X',x})$ ($\nu_1\ins\Xi$) as in (1.3.4--5). Then $\M_x$ can be determined by these $\D_{X',x}$-modules together with the morphisms $z_1$, $\dd_{z_1}$, $z_1\dd_{z_1}$ between them, where the first two morphisms are between $\M_{x,z_1}^{0}$ and $\M_{x,z_1}^{-1}$, and the last one is an endomorphism of $\M_{x,z_1}^{(\nu_1)}$ ($\nu_1\ins\Xi$).
\ms
The following is a corollary of \cite[Prop.\,II.5.4]{De}.
\msn
{\bf Proposition\hs1.3b.} {\it Let $X$ be a smooth complex variety with $Y\sst X$ a divisor with normal crossings. For a $\C$-local system $L$ on $(X\stm Y)\an$, there is uniquely a left $\D_X$-module $\M_X(L)$ such that the associated $\D_{X\an}$-module $\M_X(L)\an$ is regular holonomic of normal crossing type and we have
$$\aligned\Hc^0_{[X|Y]}\M_X(L)&=\M_X(L),\\ \DRt_X\bl(\M_X(L)\an\br)|_{(X{\setminus}Y)\an}&=L[d_X],\endaligned$$
which gives an exact functor from the category of $\C$-local systems on $(X\stm Y)\an$ to $\Msf_{\rm hol}(\D_X)$. Here $\DRt$ means that the de Rham complex is viewed as a complex of $\C$-modules on $X\an$, and we have $\M_X(L)\an\eq\M_X(L)$, $X\an\eq X$ in the analytic case. Moreover, in the algebraic case $\M_X(L)$ is regular at infinity, that is, it is isomorphic to the restriction of $\M_{\Xb}(L)$ to $X$ for any smooth compactification $\Xb$ of $X$ such that $\Xb\stm(X\stm Y)$ is a divisor with normal crossings on $\Xb$.}
\ms
(The above $\M_X(L)$ is called the meromorphic Deligne extension of a local system $L$.)
\msn
{\bf Remark\hs1.3e} ({\it Deligne extensions\1}). We first consider the {\it analytic case.} Let $z_1,\dots,z_n$ be local coordinates around $x\ins X$ as in the beginning of {\bf 1.3}. These are defined on an open subset $V\sst X$, which is identified with a polydisk $\De_{\1r}^{d_X}$ with $r>1$, where $\De\1_r$ is an open disk of radius $r$. Let $e_1,\dots,e_r$ be a basis of multivalued sections of $L$ on $V':=V\stm Y$, that is, a basis of $\Ga(\Vt',\rho^*L|_{V'})$ with $\rho:\Vt'\to V'$ a universal covering. This is identified with a basis of $L_{x_0}$, choosing a base point $x_0\ins V'$ and also a lifting of it to $\Vt'$. We may assume that $x_0$ corresponds to $(1,\dots,1)\ins\De_{\1r}^{d_X}$ via the coordinates $z_1,\dots,z_{d_X}$.
\sk
Let $T_j$ be the monodromy of $L_{x_0}$ around $\{z_j\eq0\}$ ($j\in[1,d_X]$). These commute with each other. We have $\log T_j\in{\rm End}(L_{x_0})$ such that the real part of any eigenvalue is contained in $I\subset\R$, where $I$ is either $[0,1)$ or $(-1,0]$. For a multivalued section $u\ins\Ga(\Vt',\rho^*L|_{V'})$, define
$$\ga(u):=\exp\bl(\msum_{j=1}^{d_X}\,(\log z_j)(-\log T_j)/2\pi i\br)u.
\leqno(1.3.7)$$
This determines a single valued section of $(j_{V'})_*(\OO_{V'}\sot_{\C}L|_{V'})$ with $j_{V'}:V'\into V$ the inclusion. The action of $z_j\dd_{z_j}$ corresponds to $-(2\pi i)^{-1}\log T_j$ on $\Ga(\Vt',\rho^*L|_{V'})$ by $\ga$, that is,
$$z_j\dd_{z_j}\ga(u)=\ga\bl(-(2\pi i)^{-1}(\log T_j)u\br).
\leqno(1.3.8)$$
\sk
In the case $I\eq[0,1)$, $\ga(u)$ defines a section of the {\it canonical Deligne extension\1} (see \cite[Prop.\,II.5.4]{De}), which is denoted by $\Lc_X$ in these notes. The $\ga(e_k)$ ($k\in[1,r]$) give local free generators of it. Moreover $\M_X(L)$ in Prop.\,1.3b is locally generated by the $\ga(e_k)$ ($k\in[1,r]$) over $\OX(*Y)$, and coincides with $\Lc_X(*Y)$, the localization of $\Lc_X$ along $Y$. Note that $\Lc_X$, $\M_X(L)$ are independent of the choice of local coordinates.
\sk
The above assertion is proved using (strict) Nilsson class functions, see \cite{De}. Indeed, any local section of the meromorphic (or canonical) Deligne extension is expressed as $\msum_{k=1}^r\,h_k\1e_k$ with $h_k$ (strict) Nilsson class functions. Here {\it Nilsson class functions\1} are generated over $\OX$ locally by $\mprod_{i=1}^s\,z_i^{\al_i}(\log z_i)^{p_i}$ for $\al_i\ins\C$, $p_i\ins\N$ if $Y$ is locally defined by $\mprod_{i=1}^s\1z_i\eq0$, where the condition ${\rm Re}\,\al_i\gess 0$ is added for {\it strict\1} Nilsson class functions. (This is independent of the choice of coordinates $z_1,\dots,z_{d_X}$ with $Y\eq\bl\{\mprod_{i=1}^s\1z_i\eq0\br\}$ locally.)
\sk
In the case $I\eq(-1,0]$, each $\ga(e_k)$ belongs to $\M_{x,z}^{(\nu)}$ for some $\nu\ins\Xi'\1^{d_X}$, assuming that $e_k$ is an eigenvector of the semisimple part of $\log T_j$ for any $j\in[1,d_X]$. (This is allowed, since the $\log T_j$ commute with each other.)
\sk
We now consider the {\it algebraic case.} Take a smooth compactification $\Xb$ of $X$ such that $\Xb\stm(X\stm Y)$ is a divisor with normal crossings on $\Xb$, and apply GAGA to $\Lc_{\Xb\an}$ in order to get the locally free $\OO_{\!\Xb}$-module $\Lc_{\Xb}$. We then get $\M_{\Xb}(L)$ by localizing $\Lc_{\Xb}$ along $\Xb\stm(X\stm Y)$, and $\M_X(L)$ is obtained by restricting this to $X$. We see that this is independent of the compactification $\Xb$. Indeed, if we have the Cartesian diagram
$$\begin{array}{ccccc}X&\buildrel{i_2}\over\into&\Xb_2\\ |\!|&\raise4mm\h{}\raise-2mm\h{}&\,\downarrow\!{\scriptstyle\pi}\\X&\buildrel{i_1}\over\into&\Xb_1\end{array}$$
with $\pi$ proper, then we have the isomorphisms
$$i_1^{-1}\M_{\Xb_1}(L)=i_1^{-1}\pi_*\uD\M_{\Xb_2}(L)={\rm id}_*\uD i_2^{-1}\M_{\Xb_2}(L)=i_2^{-1}\M_{\Xb_2}(L),$$
where the first isomorphism follows from Prop.\,1.3c below using Prop.\,1.2 and also {\bf A.6} below. Here we have $\Hc^j\pi_*\uD\M_{\Xb_2}(L)\eq0$ ($j\ne 0$), since this holds on $X$ and
$$\Hc^0_{[\Xb| D]}\Hc^j\pi_*\uD\M_{\Xb_2}(L)=\Hc^j\pi_*\uD\M_{\Xb_2}(L)\q(j\ins\Z),$$
with $D:=\Xb\stm(X\stm Y)$. (Recall that the localization $\Hc^0_{[\Xb| D]}$ along the divisor $D$ is an exact functor.)
\msn
{\bf Proposition\hs1.3c.} {\it Let $X,Y,L,\M_X(L)$ be as in {\rm Prop.\,1.3b}. Let $U'\sst X$ be an open subset such that $U'\,{\supset}\,U\,{:=\,}X\stm Y$. Assume the intersection of $U'$ with any irreducible component of $Y$ is non-empty. Let $\M'$ be a coherent $\D_X$-module such that $\Hc^0_{[X|Y]}\M'=\M'$ and $\M'|_{U'}\cong\M_X(L)|_{U'}$ as $\D_{U'}$-modules. Then the last isomorphism on $U'$ is extended uniquely to an isomorphism of $\D_X$-modules $\M'\cong\M_X(L)$.}
\msn
{\it Proof.} Using GAGA and a compactification, the assertion is reduced to the analytic case, and we may assume that $X$ is complex manifold. The isomorphism $\M'\cong\M_X(L)$ extending the one on $U'$ is unique (considering the image of the difference of two isomorphisms, which is supported on $Y$, hence vanishes). So the assertion is local. Using the uniqueness, we may assume that $U'$ is the {\it maximal\1} open subset satisfying the conditions in Prop.\,1.3c, that is,
$$\M'|_{U'}\cong\M_X(L)|_{U'}.$$
We will induce a contradiction assuming $U'\ne X$. (The argument is essentially equivalent to the one by induction on strata.)
\sk
Since the assertion is local, we may assume that there is a coherent $\OX$-submodule $\F\sst\M'$ generating $\M'$ over $\D_X$ (shrinking $X$). We may assume that $\F|_{U'}$ is contained in the {\it canonical extension\1} $\Lc_X$ in Rem.\,1.3e, replacing $\F$ with $\F(-mY)$ for $m\gg 0$ locally on $X$. (This inclusion can be reduced to the one on a sufficiently small open neighborhood of a general point of the intersection of $U'$ with each irreducible component of $Y$ using Rem.\,1.3f below.)
\sk
For any local generators of $\F$, their restrictions to $U'$ then determine local sections of $\Lc_X|_{U'}$, and these can be extended over an open subset $U''\sst X$ which is strictly larger than $U'$, using a Hartogs extension theorem \cite{Har} (inductively if necessary), see Rem.\,1.3f below. (Indeed, in the case $Y_{\rm sm}\sst U'$, where $Y_{\rm sm}$ denotes the smooth part of $Y$, we have $\dim X\stm U'\less d_X{-}2$, hence we get locally $U''\eq X$ by Hartogs. In the other case, we have $Y_{\rm sm}\stm U'\ne\emptyset$, and $U''$ can contain some part of it again by Hartogs.)
\sk
Shrinking $X$ if necessary, we may assume that $X\eq U''$ in order to induce a contradiction.
We then get a morphism of $\OX$-modules
$$\phi:\F\tos\M_X(L).$$
We show that this induces a morphism of $\D_X$-modules
$$\phi'':\M'\to\M_X(L),$$
using a presentation (locally on $X$)
$$\D_X\sot_{\OX}\F''\buildrel{\rho}\over\to\D_X\sot_{\OX}\F'\to\M'\to 0.$$
Here $\F'$ is a free sheaf of finite length such that the images of its free generators in $\M'$ are the above local generators of $\F$. Observe that $\phi$ induces a morphism of $\D_X$-modules
$$\phi':\D_X\sot_{\OX}\F'\to\M_X(L),$$
and its composition with $\rho$ vanishes on $U'$, since $\phi|_{U'}$ is extended to an isomorphism of $\D_{U'}$-modules $\M'|_{U'}\simto\M_X(L)|_{U'}$. This implies that the image of $\phi'\ssc\rho$ is supported on $X\stm U'\sst Y$, and it vanishes, since $\Hc^0_{[X|Y]}\M_X(L)\eq\M_X(L)$. By the above exact sequence, we then get a morphism of $\D_X$-modules
$$\phi'':\M'\to\M_X(L).$$
This is an isomorphism by Rem.\,1.2a, since $\Hc^0_{[X|Y]}\M'\eq\M'$. We thus get a contradiction. This finishes the proof of Prop.\,1.3c.
\msn
{\bf Remark\hs1.3f} ({\it Hartogs extension theorem\1}). A {\it Hartogs extension theorem\1} asserts that any holomorphic function $f$ on
$$V:=(\De\stm K){\times}U\cup\De{\times}A,$$
can be extended to a holomorphic function on $\De{\times}U$, see \cite{Har}. Here $\De$ is a unit open disk with $K\sst\De$ any compact subset, and $U$ is a connected complex manifold with $A\sst U$ any {\it non-empty open\1} subset. (A typical case is that $K\eq\{0\}$ and $U\stm A$ is a divisor, where $\De{\times}U\stm V$ has codimension 2. However, it is not necessary to assume that $U\stm A$ is an analytic subset.) This assertion can be verified by using the {\it Cauchy integral}
$$F(z,w):=\frac{1}{2\pi i}\int_{\dd\De_r}\frac{f(\zeta,w)}{\zeta{-}z}\,\ddd\zeta\q\q\bl((z,w)\ins\De\1_r{\times}U\br),$$
for any $r\ins(0,1)$ with $\De\1_r\,{\supset}\,K$ (where $\De\1_r$ is an open disk of radius $r$). Since the boundary $\dd\De\1_r$ is {\it compact,} it follows from the general theory of Lebesgue integration that $F(z,w)$ is a continuous function on $\De\1_r{\times}U$ and satisfies the {\it Cauchy-Riemann equations,} hence it is holomorphic. It coincides with $f$ on $\De\1_r{\times}A$ by the {\it Cauchy integral formula,} hence on $V\cap\De\1_r{\times}U$ by {\it analytic continuation.}
\ms
The following is essentially due to \cite[Lem.\,17]{AH} in the $L$ constant case (see also the proof of \cite[Thm.\,2]{Gr1}) and \cite[Prop.\,II.6.8]{De} in general. We note here a short proof for the convenience of the reader.
\msn
{\bf Proposition\hs1.3d.} {\it Let $Y$ be a divisor with normal crossings on a complex manifold $X$. For a $\C$-local system $L$ on $U:=X\stm Y$, let $\M_X(L)$ be the meromorphic Deligne extension of $L$, see {\rm Prop.\,1.3b}. Let $j_U:U\into X$ be the canonical inclusion. Then we have the quasi-isomorphisms}
$$\DRt_X\bl(\M_X(L)\br)\simto\RR(j_U)_*j_U^*\DRt_X\bl(\M_X(L)\br)\simot\RR(j_U)_*L[d_X].
\leqno(1.3.9)$$
\msn
{\it Proof.} It is enough to show the first isomorphism, since the second is clear. Since the assertion is local, we may assume that $X$ is a polydisk $\De^{d_X}$. There is a filtration $G$ on $L$ such that the graded pieces are local systems of rank 1, since the local monodromies commute with each other. Using the distinguished triangles associated to the short exact sequences $0\tos G_{p-1}\tos G_p\tos\Gr^G_p\tos 1$ inductively, the assertion is reduced to the case $L$ has rank 1.
\sk
Let $T_k$ be the monodromy around $\{z_k\eq0\}$, which is identified with a complex number $\la_k$ ($k\in[1,n]$). Then $H\ub(U,L)$ is calculated by the Koszul complex for the multiplications by $\la_1,\dots,\la_n$ on $\C$. (This well-known in the case $n\eq1$. The general case can be reduced to this case inductively by using the projection $\De^{d_X}\tos\De^{d_X-1}$.)
\sk
Let $\Lc_X$ be the canonical extension of $\OO_U\sot_{\C}L$ in \cite[Prop.\,II.5.4]{De} such that the real parts of eigenvalues of residues of connection (which are $\al_k\eq{-}\log\la_k/2\pi i$ by (1.3.8)) are contained in $[0,1)$, see Rem.\,1.3e. We have the logarithmic de Rham complex denoted by $\DR_X^{\log}(\Lc_X)$. This is a subcomplex of $\DR_X\bl(\M_X(L)\br)$. It is easy to show the quasi-isomorphism
$$\DR_X^{\log}(\Lc_X)\simto\DR_X\bl(\M_X(L)\br),
\leqno(1.3.10)$$
since these are the Koszul complexes associated with the actions of $z_k\dd_{z_k}$ ($k\less n$) and $\dd_{z_k}$ ($k\,{>}\,n$) on $\Lc_X$ and $\M_X(L)$ respectively, see also \cite{De}, \cite[3.11]{mhm}. Indeed, we have the filtrations $F^{(i)}$ ($i\in[1,n]$) on $\M_X(L)_x$ defined by
$$\aligned&F^{(i)}_p\M_X(L)_x:=\M_X(L)_x\cap\mprod_{{\rm Re}\,\nu_i\ges-p}\,\M_X(L)_{x,z}^{(\nu)},\\&\h{so that}\q F^{(1)}_{p_1}\cdots F^{(n)}_{p_n}\M_X(L)_x\eq z_1^{-p_1}\cdots z_n^{-p_n}\Lc_{X,x}.\endaligned$$
\sk
The assertion (1.3.9) is furthermore reduced to the case $n\eq d_X$, using the associated double complex structure of the Koszul complex together with the Poincar\'e lemma for the $\dd_{z_k}$ ($k\,{>}\,n$). We also see that the stalk of the complexes at $0\ins\De^{d_X}$ are acyclic if $\al_k\ne 0$, that is, $\la_k\ne 1$, for some $k\in[1,n]$. The assertion is thus reduced to the case $\la_k\eq 1$ for any $k\in[1,n]$, that is, $L\eq\C_U$, and $n\eq d_X$.
\sk
It is then enough to prove the canonical isomorphisms
$$\Hc^j\DR^{\log}_X(\OX)_0\simto\rlap{\raise-10pt\h{$\,\,\,\scriptstyle\de$}}\indlim\, H^j(U_{\de},\C)\q(j\ins\Z).
\leqno(1.3.11)$$
Here $U_{\de}:=(\De_{\de}^*)^n$ with $\De_{\de}^*$ the punctured disk of radius $\de\,{>}\,0$, and the above inductive system is constant. The left-hand side of (1.3.11) is the cohomology of the Koszul complex associated with the actions of $z_k\dd_{z_k}$ on $\C\{z_1,\dots,z_n\}$. Since each monomial is stable by these actions, we can verify that this Koszul complex is quasi-isomorphic to the Koszul subcomplex associated with the 0 actions on $\C$. (Here we have to examine some convergence. This is closely related to the Poincar\'e lemma for convergent power series.) So both sides of (1.3.11) have the same dimension. It is then sufficient to show the surjectivity. However, this can be verified using the integration of logarithmic $k$-forms $\ddd z_{i_1}/z_{i_1}{\wedge}\cdots{\wedge}\ddd z_{i_k}/z_{i_k}$ on topological $k$-cycles $\ga_{j_1}{\times}\cdots{\times}\ga_{j_k}\sst U_{\de}$, which does not vanish if and only if $\{i_1,\dots,i_k\}\eq\{j_1,\dots,j_k\}$, where $\ga_j:=\{|z_j|\eq\de/2\}\sst\De^*_{\de}$. So (1.3.11) follows. This finishes the proof of Prop.\,1.3d.
\ms
Slightly generalizing the above argument, we can show the following.
\msn
{\bf Proposition\hs1.3e.} {\it Let $f:X\to S$ be a proper smooth morphism of complex manifolds. Let $Y\sst X$ be a divisor with relative normal crossings, that is, $Y$ is a divisor with normal crossings and any intersection of local irreducible components of $Y$ is smooth over $S$. Let $L$ be a $\C$-local system on $U:=X\stm Y$ with $\M_X(L)$ as in {\rm Prop.\,1.3b}. Let $f_U$ be the restriction of $f$ to $U$. Then the $\Hc^jf_*\uD\M_X(L)$ are locally free $\OO_S$-modules of finite rank such that}
$$\DR_S\bl(\Hc^jf_*\uD\M_X(L)\br)[-d_S]=R^{j+d_X-d_S}(f_U)_*L.
\leqno(1.3.12)$$
\msn
{\it Proof.} Let $t_1,\dots,t_{d_S}$ be local coordinates of $S$. Set $f_i:=f^*t_i$. Let $i_f:X\into X{\times}S$ be the graph embedding of $f$ with ${\rm pr}_2:X{\times}S\tos S$ the second projection. We have the isomorphisms
$$\aligned f\uD_*\M_X(L)&=\RR({\rm pr}_2)_*\DR_{X\times S/S}\bl((i_f)\uD_*\M_X(L)\br),\\(i_f)\uD_*\M_X(L)&=\M_X(L)[\dd_{t_1},\dots,\dd_{t_{d_S}}]\de(t_1{-}f_1)\cdots\de(t_{d_S}{-}f_{d_S}),\endaligned$$
with $\dd_{z_j}\de(t_i{-}f_i)=-(\dd_{z_j}f_i)\dd_{t_i}\de(t_i{-}f_i)$ (by an argument similar to \cite[1.1]{rp}).
\sk
Let $F$ be the filtration on the complex
$$\DR_{X\times S/S}\bl((i_f)\uD_*\M_X(L)\br),$$
by the total order of the $\dd_{t_i}$ shifted by the complex degree. Its graded quotients are then isomorphic to the complex
$$\bl(\M_X(L)[\dd_{t_1},\dots,\dd_{t_{d_S}}]\sot_{\OX}\Om_X\ub[d_X],\,-\msum_{i=1}^{d_S}\,\dd_{t_i}\ddd f_i\wedge\br).$$
Take local coordinates $z_1,\dots,z_{d_X}$ of $X$ with $z_i\eq f_i$ for $i\in[1,d_S]$. We can verify that the above filtered complex is filtered quasi-isomorphic to its subcomplex endowed with the induced filtration and such that their components are the intersections of the kernels of the actions of $\ddd z_i\wedge$ ($i\in[1,d_S]$) on $\M_X(L)\sot_{\OX}\Om_X\ub[d_X]$. (Here we use also the filtration by the order of each $\dd_{t_p}$ shifted by the complex degree, and proceed by decreasing induction on $p$ considering the intersections of the kernels of $\ddd z_i\1\wedge$ for $i\,{>}\,p$ together with the coefficient polynomial rings $\C[\dd_{t_1},\dots,\dd_{t_p}]$.)
\sk
Taking the division by $\ddd z_1{\wedge}\cdots{\wedge}\ddd z_{d_S}$, we can then see that the above subcomplex is isomorphic as a complex of $f^{-1}\OO_S$-modules to the relative de Rham complex
$$\DR_{X/S}\bl(\M_X(L)\br),$$
whose \h{$k$\1th} component is
$$\M_X(L)\sot_{\OX}\Om_{X/S}^{k+d_X-d_S}.$$
\sk
By an argument similar to the proof (1.3.10) (that is, with parameter), we can verify that the natural inclusion induces the quasi-isomorphism
$$\DR^{\log}_{X/S}(\Lc_X)\simto\DR_{X/S}\bl(\M_X(L)\br),$$
where $\DR^{\log}_{X/S}$ is the relative logarithmic de Rham complex with $\Lc_X$ as in Rem.\,1.3e.
\sk
We have thus proved the following isomorphism in the derived category of complexes of $f^{-1}\OO_S$-modules:
$$f\uD_*\M_X(L)\cong\DR^{\log}_{X/S}(\Lc_X).$$
Applying the Grauert coherence theorem (see for instance \cite{Dou}), we then see that the $\Hc^jf\uD_*\M_X(L)$ ($j\ins\Z$) are coherent $\OO_S$-modules, hence locally free. Indeed, it is well-known that a $\D_S$-module is locally free over $\OO_S$ if it is coherent over it. (This easily follows from {\bf A.7} below by reducing to the curve case.) The remaining assertion then follows from Prop.\,1.3d using the commutativity of the de Rham functor with the direct image functors. (This commutativity can be verified easily in the projection case. The local case is also easily proved if we choose local coordinates; however, the independence of the isomorphism under the choice of local coordinates does not seem trivial. Here we can use the theory of induced $\D$-modules associated to differential complexes together with the forgetful functor, see for instance \cite[3.1]{ind}, etc.) This finishes the proof of Prop.\,1.3e.
\msn
{\bf 1.4.~Regular meromorphic connections.} Let $Z$ be a complex variety with $D\sst Z$ a locally principal divisor such that $Z\stm D$ is smooth. Let $\pi:(\Zt,E)\tos(Z,D)$ be an embedded resolution such that $E:=\pi^{-1}(D)\sst\Zt$ is a divisor with normal crossings and $\pi$ induces an isomorphism over $Z\stm D$. (We may assume that $\pi$ is projective.)
\sk
We say that $\M$ is a {\it meromorphic connection with poles along\1} $D$ if $\M$ is a coherent $\OZ(*D)$-module endowed with a sheaf morphism
$$\nabla:\M\to\Om_Z^1(*D){\otimes}_{\OZ}\M,
\leqno(1.4.1)$$
(which is literally a meromorphic connection), and its restriction to $Z\stm D$ is an {\it integrable connection\1} in the usual sense. Here $\OZ(*D)$ is the localization of $\OZ$ by a local defining function of $D$. Note that coherence over $\OZ(*D)$ is equivalent to that it is quasi-coherent over $\OZ$, and is locally finitely generated over $\OZ(*D)$.
\sk
If $\M$ is a meromorphic connection on $Z$ with poles along $D$, then so is $\pi^*\!\!\M$ on $\Zt$ with poles along $E$. (Here $\pi_{\OO}^*$ is denoted by $\pi^*$ to simplify the notation.) Indeed, we have the induced connection
$$\pi^*\nabla:\pi^*\!\!\M\to\Om_{\Zt}^1(*E){\otimes}_{\OO_{\!\Zt}}\pi^*\!\!\M,
\leqno(1.4.2)$$
since
$$\pi^*\bl(\Om_Z^1(*D){\otimes}_{\OZ}\M\br)=\pi^*\Om_Z^1(*D){\otimes}_{\OO_{\!\Zt}}\pi^*\!\!\M.
\leqno(1.4.3)$$
Here the canonical morphism $\pi^*\Om_Z^1\to\Om_{\Zt}^1$ induces the isomorphism
$$\pi^*\Om_Z^1(*D)\simto\Om_{\Zt}^1(*E).
\leqno(1.4.4)$$
(considering its kernel and cokernel which are supported in $E$).
\sk
In the $Z$ smooth case (for instance, on $\Zt\1$), meromorphic connections are equivalent to left $\D_Z(*D)$-modules which are coherent over $\OZ(*D)$. Here $\D_Z(*D)$ is the ring locally generated by $\D_Z$ and $g^{-1}$ in $(j_{Z\setminus D})_*j_{Z\setminus D}^{-1}\D_Z$ with $g$ a local defining function of $D\sst Z$.
\sk
In the normal crossing divisor case (that is, on $\Zt$), a meromorphic connection $\M$, or its corresponding left $\D_Z(*D)$-module which is coherent over $\OZ(*D)$, is called {\it regular\1} if it is isomorphic to the meromorphic Deligne extension of a local system as in Prop.\,1.3b. In general, a meromorphic connection $\M$ on $Z$ with poles along $D$ is called {\it regular\1} if so is $\pi^*\!\!\M$.
Set
$$\D_Z(*D):=\pi_*\D_{\Zt}(*E).
\leqno(1.4.5)$$
This is independent of the choice of a desingularization. Indeed, in case $Z$ is smooth, $\pi_*F_p\D_{\Zt}$ is a coherent subsheaf of $F_p\D_Z$ by the Hartogs extension theorem, and the two sheaves coincide after the localization along $D$ for each $p\ins\Z$. We have
$$\D_Z(*D)\supset\OZ(*D).
\leqno(1.4.6)$$
\sk
Note that (1.4.4) implies the isomorphism
$$\Om_Z^1(*D)\simto\pi_*\pi^*\Om_Z^1(*D)\simto\pi_*\Om_{\Zt}^1(*E).
\leqno(1.4.7)$$
We have also the isomorphism
$$\Theta_Z(*D)=\pi_*\Theta_{\Zt}(*E),
\leqno(1.4.8)$$
localizing along $D$ the canonical morphisms
$$\pi_*\Theta_{\Zt}\to\pi_*\pi^*\Theta_Z\leftarrow\Theta_Z.$$
\sk
It follows from (1.4.5), (1.4.8) that a left $\D_Z(*D)$-module which is coherent over $\OZ(*D)$ is a meromorphic connection on $Z$ with poles along $D$. Its converse follows from this assertion on $\Zt$ (which is easy and well-known) using (1.4.9) below.
\sk
A left $\D_Z(*D)$-module which is coherent over $\OZ(*D)$ is called {\it regular\1} if it is so as a meromorphic connection, that is, if its pullback to $\Zt$ as a meromorphic connection is regular or if it is isomorphic to the direct image of a regular meromorphic connection on $\Zt$ using (1.4.9--10) below.
\sk
For a meromorphic connection $\M$ on $Z$ with poles along $D$, we have the isomorphisms
$$\pi_*\pi^*\!\!\M=\M,\q\pi_*\bl(\Om_{\Zt}^1(*E){\otimes}_{\OO_{\!\Zt}}\pi^*\!\!\M\br)=\Om_Z^1(*D){\otimes}_{\OZ}\M,
\leqno(1.4.9)$$
using (1.4.7) for the last isomorphism.
\sk
For a {\it regular\1} meromorphic connection $\M'$ on $\Zt$ with poles along $E$, we have the coherence of $\pi_*\M'$ over $\OZ(*D)$ by the Grauert coherence theorem applied to the Deligne canonical extension. There are isomorphisms
$$\pi^*\pi_*\M'=\M',\q\pi^*\bl(\Om_Z^1(*D){\otimes}_{\OZ}\pi_*\M'\br)\simto\Om_{\Zt}^1(*E){\otimes}_{\OO_{\Zt}}\M'.
\leqno(1.4.10)$$
\sk
Combined with \cite[Prop.\,II.5.4]{De}, the above arguments imply the following.
\msn
{\bf Proposition\hs1.4a.} {\it There are equivalences between the abelian categories with objects as follows}\,$:$
\skn
\cond[(a)] {\it $\C$-local systems on $Z\an\stm D\an$,}
\skn
\cond[(b)] {\it Vector bundles with $($regular$\,)$ integrable connection on $Z\stm D$,}
\skn
\cond[(c)] {\it Regular meromorphic connections on $\Zt$ with poles along $E$,}
\skn
\cond[(d)] {\it Regular left $\D_{\Zt}(*E)$-modules which are coherent over $\OO_{\Zt}(*E)$,}
\skn
\cond[(e)] {\it Regular meromorphic connections on $Z$ with poles along $D$,}
\skn
\cond[(f)] {\it Regular left $\D_Z(*D)$-modules which are coherent over $\OZ(*D)$.}
\ms
Note that the ``regularity" (at infinity) in (b) is needed only for the {\it algebraic\1} case. In (a) $Z\an\stm D\an=Z\stm D$ in the analytic case.
\sk
We have moreover the following.
\msn
{\bf Proposition\hs1.4b} ({\it Embedded case\1}).~{\rm (1)} {\it Assume $Z$ is a reduced closed subvariety of a smooth variety $X$, and there is a divisor $Y\sst X$ such that $D\eq Y\caps Z$ set-theoretically. The category for {\rm(e)} in {\rm Prop.\,1.4a} with regularity condition forgotten is equivalent to the category with objects as follows}\,$:$
\skn
\cond[{\rm(g)}] {\it Coherent left $\D_X(*Y)$-modules whose characteristic varieties in the cotangent bundle of $X\stm Y$ are contained in the conormal bundle of $Z\stm D$ in $X\stm Y$.}
\sk
{\rm (2)} {\it An object of {\rm (g)} corresponds to an object of {\rm(e)} in {\rm Prop.\,1.4a} $($with regularity condition {\rm not} forgotten$)$ if and only if it is isomorphic to the direct image as a $\D$-module of a regular left $\D_{\Zt}(*E)$-module which is coherent over $\OO_{\!\Zt}(*E)$.}
\msn
{\it Proof.} To show the assertion (1), it is better to use right $\D$-modules rather than left $\D$-modules. Let $\M^r$ be a right $\D_Z(*D)$-module which is coherent over $\OZ(*D)$. Then we can define the {\it direct image\1} of $\M^r$ by
$$(i_Z)_*^{\D}\M^r:=\M^r{\otimes}_{\D_Z(*D)}\1\D_{Z\to X}(*Y),
\leqno(1.4.11)$$
setting
\vskip-8mm
$$\D_{Z\to X}(*Y):=\pi_*\D_{\Zt\to X}(*Y),$$
which is a left $\D_Z(*D)$ and right $i_Z^{-1}\D_X(*Y)$-bimodule. The sheaf-theoretic direct image by the inclusion $i_Z:Z\into X$ is omitted to simplify the notation.
\sk
Let $\Nc^r$ be a coherent right $\D_X(*Y)$-module such that its characteristic variety (defined in the cotangent bundle of $X\stm Y$) is contained in the conormal bundle of $Z\stm D$ in $X\stm Y$. Set
$$\Hc^0(i_Z)_{\D}^!\Nc^r:=(\Nc^r)^{\I_Z}=\Hom_{\OX(*Y)}(\OZ(*D),\Nc^r).
\leqno(1.4.12)$$
Here $\I_Z\sst\OX$ is the reduced ideal sheaf of $Z$ so that $\OZ\eq\OX/\!\1\I_Z$, and $(\Nc^r)^{\I_Z}$ denotes the subsheaf annihilated by $\I_Z$. (This functor coincides with the pullback functor $\Hc^0(i_Z)_{\D}^!$ defined in {\bf A.4} on $X\stm Y$.)
\sk
Applying Rem.\,A.7b, these two operations are inverse of each other. Indeed, there are canonical morphisms
$$\M^r\to\Hc^0(i_Z)_{\D}^!(i_Z)_*^{\D}\M^r,\q(i_Z)_*^{\D}\1\Hc^0(i_Z)_{\D}^!\Nc^r\to\Nc^r.
\leqno(1.4.13)$$
inducing isomorphisms on $Z\stm D$ and $X\stm Y$ respectively (hence everywhere).
\sk
Here we have to show the finiteness of $(\Nc^r)^{\I_Z}$ over $\OZ(*D)$ assuming the coherence of $\Nc^r$ over $\D_X(*Y)$ together with the condition on the characteristic variety. Take locally a coherent filtration $F$ of a coherent $\D_X$-submodule $\,{}'\!\!\!\Nc^r\sst\Nc^r$ whose localization along $Y$ is $\Nc^r$. The supports of $(\,\,'\!\!\!\Nc^r)^{\I_Z}/(F_p\,\,{}'\!\!\!\Nc^r)^{\I_Z}$ for $p\ins\N$ form a decreasing sequence of closed subvarieties of $X$, and are contained in $Y$ for $p\gg 0$ locally on $X$. So the assertion~(1) follows.
\sk
To show the assertion (2), let $\M'$ be a regular meromorphic connection on $\Zt$ with poles along $E$. We have the canonical isomorphism
$$\aligned\om_Z(*D){\otimes}_{\OZ}\pi_*\M'\simto{}&\pi_*\bl(\om_{\Zt}(*E){\otimes}_{\OO_{\!\Zt}}\M'\br)\\=\,\,{}&\pi_*\M'{}^r,\endaligned
\leqno(1.4.14)$$
setting
\vskip-7mm
$$\om_Z(*D):=\pi_*\om_{\Zt}(*E),\q\M'{}^r:=\om_{\Zt}(*E){\otimes}_{\OO_{\!\Zt}}\M'.$$
Note that the right-hand side of (1.4.14) is by definition the right $\D_Z(*D)$-module $\M^r$ corresponding to the meromorphic connection $\M:=\pi_*\M'$. (It is unclear whether $\om_Z(*D)$ is locally free over $\OZ(*D)$, but it is nevertheless {\it invertible.})
\sk
We have moreover the canonical isomorphism
$$\aligned\pi_*\M'{}^r{\otimes}_{\D_Z(*D)}\1\D_{Z\to X}(*Y)\simto{}&\pi_*\bl(\M'{}^r{\otimes}_{\D_{\Zt}(*E)}\1\D_{\Zt\to X}(*Y)\br)\\=\,\,{}&(i_Z\ssc\pi)_*^{\D}\M'{}^r.\endaligned
\leqno(1.4.15)$$
So the assertion (2) follows from the isomorphisms in (1.4.9--10) and (1.4.13). This finishes the proof of Prop.\,1.4b.
\msn
{\bf 1.5.~Relative logarithmic de Rham complexes over a disk.} Let $f:X\to\De$ be a nonconstant morphism from a connected complex manifold $X$ onto a disk $\De$. Let $Y\sst X$ be a divisor with normal crossings containing $X_0:=f^{-1}(0)\sst X$. (Note that $X_0$ is also a divisor with normal crossings, since it is a divisor contained in $Y$.) For a $\C$-local system $L$ on $U:=X\stm Y$, we denote by $\M_X(L)$ the meromorphic Deligne extension of $L$, see Prop.\,1.3b. Let $i_f:X\into X{\times}\De$ be the graph embedding with $t$ the coordinate of $\De$. The direct image of $\M_X(L)$ by $i_f$ as a left $\D$-module can be described as
$$(i_f)_*\uD\M_X(L)=\M_X(L)[\dd_t]\1\de(f{-}t),$$
with actions of $\dd_{z_i},t$ on $\de(f{-}t)$ defined by
$$\dd_{z_i}\de(f{-}t)\eq{-}(\dd_{z_i}f)\1\dd_t\1\de(f{-}t),\q t\1\de(f{-}t)\eq f\1\de(f{-}t).$$
Here the actions of $\OX,\dd_t$ are given naturally, see for instance \cite[1.1]{rp}.
\sk
The direct image of $\M_X(L)$ by $f$ is then expressed as
$$f_*\uD\M_X(L)=\RR f_*\DR_X\bl(\M_X(L)[\dd_t]\1\de(f{-}t)\br).$$
\sk
We denote by $\Lc_X$ the canonical extension of $\OO_U\sot_{\C}L$ in \cite[Prop.\,II.5.4]{De} with real part of any eigenvalue of the residues of connection contained in $[0,1)$. We have the isomorphism $\M_X(L)\eq\Lc_X(*Y)$, which implies that
$$\DR_X^{\log}\bl(\Lc_X(*Y)[\dd_tt]\1\de(f{-}t)\br)=\DR_X\bl(\M_X(L)[\dd_t]\1\de(f{-}t)\br),$$
where $\DR_X^{\log}$ is defined by using the logarithmic complex $\Om_X\ub\lYr$ (instead of $\Om_X\ub$), see \cite{De}. Note that the differential of the complex is given by
$$\nabla\mi\dd_t\1\ddd f\wedge=\nabla\mi\dd_tt\1\tfrac{\ddd f}{f}\wedge,
\leqno(1.5.1)$$
where $\nabla$ is associated to the logarithmic connection on $\Lc_X(*Y)$.
\sk
Take local coordinates $z_1,\dots,z_{d_X}$ such that $f\eq\mprod_{i=1}^r\,z_i^{m_i}$ and $Y$ is locally defined by $\mprod_{i=1}^n\,z_i\eq0$, where $1\less r\less n\less d_X$. In the notation of {\bf 1.3}, set
$$V:=\mopl_{\nu\in\Xi'{}^n}\,\M_X(L)_{x,z}^{(\nu)},$$
where $\Xi'{}^n$ is identified with a subset of $\Xi'{}^{d_X}$ by adding $(0,\dots,0)\ins\Xi'{}^{d_X-n}$ if $n\,{<}\,d_X$. We denote by $R_i$ the action of $z_i\dd_{z_i}$ on $V$ if $i\less n$, and 0 otherwise.
\sk
We can then locally identify $\DR_X^{\log}\bl(\Lc_X(*Y)[\dd_tt]\1\de(f{-}t)\br)$ with the shifted Koszul complex associated with the actions on $\OX(*Y)\sot_{\C}V[\dd_t]$ defined by
$$\xi_i:=\begin{cases}z_i\dd_{z_i}{+}\1 R_i{-}\1 m_if\1\dd_t&(i\less n),\\ \dd_{z_i}&(i\,{>}\, n),\end{cases}$$
using the local generators $\ddd z_i/z_i$ ($i\less n$), $\ddd z_i$ ($i\,{>}\,n$) of $\Om_X^1\lYr$ and their exterior products. Here $m_i\eq 0$ for $i\,{>}\,r$, and $z_i\dd_{z_i}f\eq m_if$ for any $i$.
\sk
Consider the following subcomplexes of $\DR_X^{\log}\bl(\Lc_X(*Y)[\dd_tt]\1\de(f{-}t)\br)$:
$$\aligned\K_f\ub(\Lc_X)&\,\subset\,\DR_X^{\log}\bl(\Lc_X[\dd_tt]\1\de(f{-}t)\br)\q\q\q\h{with}\\\K_f^j(\Lc_X)&:=\Lc_X\,\sot_{\OX}\,\K_f^j,\\ \K_f^j&:={\rm Ker}\bl(\tfrac{\ddd f}{f}\wedge:\Om_X^{d_X+j}\lYr\to\Om_X^{d_X+j+1}\lYr\br).\endaligned$$
(Note that $\Lc_X\,\sot_{\OX}$ is an exact functor, and commutes with the kernel of $\tfrac{\ddd f}{f}\wedge$.)
\sk
We can show that the natural inclusion
$$\iota:\K_f\ub(\Lc_X)\into\DR_X^{\log}\bl(\Lc_X[\dd_tt]\1\de(f{-}t)\br)
\leqno(1.5.2)$$
is a quasi-isomorphism. This can be verified by using the filtration by the order of $\dd_tt$ on $\DR_X^{\log}\bl(\Lc_X[\dd_tt]\1\de(f{-}t)\br)$ which is shifted by the degree of complex. Indeed, its $k$\1th graded piece is isomorphic to the truncated shifted Koszul complex
$$\sigma^{\ges-k}\bl(\Lc_X\sot_{\OX}\Om_X\ub\lYr[d_X],\tfrac{\ddd f}{f}\wedge\br)\q(k\ins\N).$$
(Recall that, for a complex $K\ub$ in general, we have $(\sigma^{\ges k}K\ub)^j\eq K^j$ if $j\gess k$, and 0 otherwise.)
The complex $\bl(\Om_X\ub\lYr,\tfrac{\ddd f}{f}\wedge\br)$ is identified with the Koszul complex for multiplications by $m_i$ on $\OX$ ($i\ins[1,d_X]$) using $\ddd z_i/z_i$ ($i\less n$), $\ddd z_i$ ($i\,{>}\,n$). Hence it is acyclic (since $r\gess 1$). This holds also after taking the tensor product with $\Lc_X$. 
\sk
Notice that $\K_f\ub(\Lc_X)$ is a subcomplex of $\DR_X^{\log}(\Lc_X)$ which is identified with the shifted Koszul complex associated with the actions on $\Lc_X$ defined by
$$\xi'_i:=\begin{cases}z_i\dd_{z_i}{+}\1 R_i&(i\less n),\\ \dd_{z_i}&(i\,{>}\, n).\end{cases}$$
\sk
Using the acyclicity of $\bl(\Om_X\ub\lYr,\tfrac{\ddd f}{f}\wedge\br)$, we get the isomorphism of complexes of $f^{-1}\OO_S$-modules
$$\K_f\ub(\Lc_X)\simot\DR_{X/\De}^{\log}(\Lc_X):=\Om_{X/\De}\ub\lYr(\Lc_X)[d_X{-}1],
\leqno(1.5.3)$$
induced by $\tfrac{\ddd f}{f}\wedge$ (since $\tfrac{\ddd f}{f}$ is a {\it closed\1} form). Here $\Om_{X/\De}^j\lYr(\Lc_X)$ is defined by
$${\rm Coker}\bl(\tfrac{\ddd f}{f}\wedge:\Lc_X\sot_{\OX}\Om_X^{j-1}\lYr\to\Lc_X\sot_{\OX}\Om_X^j\lYr\br).$$
(This is compatible with the tensor product with $\Lc_X$.)
\sk
Under the induced isomorphisms
$$\Hc^j\bl(\K_f\ub(\Lc_X)\br){}_x\simot\Hc^j\bl(\DR_{X/\De}^{\log}(\Lc_X)\br){}_x\q(x\ins X_0),$$
the action of $t\dd_t$ on the right-hand side corresponds to that of $\dd_tt$ on the left-hand side, since the isomorphism is induced by $\tfrac{\ddd f}{f}\wedge$ (instead of $\ddd f\wedge$) and $f\eq f^*t$. Note that the action of $t\dd_t$ on the right-hand side is defined by the {\it connection morphism\1} of the long exact sequence associated with the short exact sequence of complexes
$$0\to\DR_{X/\De}^{\log}(\Lc_X)[-1]\buildrel{\!\!\theta\wedge}\over\longrightarrow\DR_X^{\log}(\Lc_X)\to\DR_{X/\De}^{\log}(\Lc_X)\to 0,
\leqno(1.5.4)$$
where the {\it first injection\1} is induced by the exterior product with $\theta\eq\tfrac{\ddd f}{f}$. This short exact sequence is isomorphic to
$$0\to\K_f\ub(\Lc_X)\to\DR_X^{\log}(\Lc_X)\buildrel{\!\!\theta\wedge}\over\longrightarrow\K_f\ub(\Lc_X)[1]\to 0,$$
where the {\it last surjection\1} is induced by the exterior product with $\theta\eq\tfrac{\ddd f}{f}$. This defines the action of $\dd_tt$ on $\Hc^j\bl(\K_f\ub(\Lc_X)\br){}_x$ using the connecting morphism. By (1.5.1) this action is compatible with the action of $\dd_tt$ on the $\Hc^j\DR_X^{\log}\bl(\Lc_X[\dd_tt]\1\de(f{-}t)\br){}_x$ via the inclusion (1.5.2).
\msn
{\bf Remark\hs1.5a.} It does not seem trivial to show the compatibility of the action of $t\dd_t$ on the direct image $\D$-modules $\Hc^jf_*\uD\M_X(L)$ ($j\ins\Z$) with the action of the {\it logarithmic Gauss-Manin connection\1} $\nabla_{\!t\dd_t}$ on the {\it relative de Rham cohomology sheaves\1} $R^j\!f_*\DR_{X/\De}^{\log}(\Lc_X)$, where the latter is defined by using the connecting morphism of the above short exact sequence after applying $\RR f_*$ to it.
\msn
{\bf Remark\hs1.5b} (One-dimensional local Riemann-Hilbert correspondence). It is well-known that there are equivalences between the categories with objects as follows:
\skn
\cond[(a)] Finite free $\C\{t\}$-modules $M$ with an action of logarithmic connection $\nabla_{\!t\dd_t}$ such that the real part of any eigenvalue of the residue is contained in $[0,1)$.
\skn
\cond[(b)] $\C$-local systems $L$ on $\De^*$.
\skn
\cond[(c)] Finite dimensional $\C$-vector spaces $V$ with an automorphism $T$.
\skn
\cond[(d)] Finite dimensional $\C$-vector spaces $V$ with an endomorphism $R$ such that the real part of any eigenvalue is contained in $[0,1)$.
\ms
Indeed, the functor from (a) to (b) is induced by the restriction of the shifted de Rham functor $\DR_{\De}[-1]$. Here a finite $\C\{t\}$-module $M$ can be extended to a coherent $\OO$-module on $\De_{\de}$ for $0\,{<}\,\de\,{\ll}\,1$, and the category of local systems on $\De_{\de}^*$ is independent of $\de\,{>}\,0$.
\sk
The functor from (b) to (c) is given by the nearby cycle functor $\psi_t$ which is naturally associated with the monodromy transformation $T$.
\sk
The functor from (c) to (d) is obtained by taking an appropriate logarithm of $T$ divided by $-2\pi i$.
\sk
We have also the functor form (a) to (d) defined by
$$(M,\nabla_{\!t\dd_t})\,\mapsto\,\mopl_{{\rm Re}\,\al\in[0,1)}\,(M^{(\al)},\nabla_{\!t\dd_t}),$$
where $M^{(\al)}:={\rm Ker}\,(t\dd_t{-}\al)^m\subset M$ ($m\gg 0$). (This is a special case of Rem.\,1.3b.)
\msn
{\bf Proposition\hs1.5a.} {\it With the notation in the beginning of {\bf 1.5}, set
$$\Hc^j_{X/\De}(\Lc_X):=\Hc^j\bl(\DR_{X/\De}^{\log}(\Lc_X)\br)\q\q(j\ins\Z).$$ 
These are constructible sheaves of finite free $\C\{t\}$-modules on $X_0$ endowed with the action of $\nabla_{\!t\dd_t}$ such that the real part of any eigenvalue of the residue is contained in $[0,1)$. Moreover there are canonical isomorphism of $\C$-constructible sheaves compatible with the action of monodromy associated with the nearby cycle functors
$$\psi_t\DR_{\De}\bl(\Hc^j_{X/\De}(\Lc_X)|_{X_0}\br)[-1]\simto\Hc^{j+d_X-1}\psi_f\RR(j_U)_*L,
\leqno(1.5.5)$$
where $j_U:U\into X$ is the natural inclusion $($see {\rm Rem.\,1.3b} for $\psi_t\DR_{\De}[-1])$.}
\msn
{\it Proof.} For $x\ins X_0$, take local coordinates $z_1,\dots,z_{d_X}$ around $x$ such that $f\eq\mprod_{i=1}^r\,z_i^{m_i}$ $(m_i\gess 1$), and $Y$ is locally defined by $\mprod_{i=1}^n\,z_i\eq 0$ with $1\less r\less n\less d_X$. Since the assertion is local, we may assume
$$X\eq\mprod_{i=1}^{d_X}\,\De_{\ep_i}\q(0\,{<}\,\ep_i\,{\ll}\,1).$$
We have $X\eq X'{\times}X''$ with $X'\eq\mprod_{i=1}^r\,\De_{\ep_i}$, $X''\eq\mprod_{i=r+1}^{d_X}\,\De_{\ep_i}$. Set
$$\aligned&X_c:=\{f\eq c\}\sst X,\q X'_c:=\{f|_{X'}\eq c\}\sst X'\\&\h{so that}\q\q
X_c=X'_c{\times}X''\q\q(c\ins\C).\endaligned$$
For $c\ins\C^*$, $X'_c$ is an unramified covering of an open subset $U_{|c|,\ep}\subset\mprod_{i=2}^r\,\De_{\ep_i}$ defined by the condition
$$\mprod_{i=2}^r\,|z_i|^{m_i}>|c|/\ep_1^{m_1}.$$
For $|c|\ll\mprod_{i=1}^r\,\ep_i^{m_i}$, we see that the cohomology groups $H^j(X_c,\RR(j_U)L|_{X_c})$ are independent of $c$, $\ep_i$ ($i\ins[1,d_X])$ using the coordinate change $z_i\,{\mapsto}\,a_iz_i$ for $a_i\ins\C^*$. Here we have to show that the natural inclusions induce isomorphisms. This can be shown by using the integral curves of an appropriate real vector field $\xi$ which is a linear combination of the $\xi_i:=(x_i\dd_{x_i}\pl y_i\dd_{y_i})/2m_i$ with $C^{\infty}$ coefficients (using a partition of unity) and such that $\xi\1|f|^2=|f|^2$. Here $z_i\eq x_i\pl\sqrt{-1}\1y_i$ so that $|z_i|^2\eq x_i^2\pl y_i^2$ and $\xi_i\1|f|^2=|f|^2$. We then get the isomorphisms
$$\Hc^j\bl(\psi_f\RR(j_U)L\br){}_x=H^j\bl(X_c,\RR(j_U)L|_{X_c}\br)\q(j\ins\Z).$$
\sk
For $0\,{<}\,\de\ll\mprod_{i=1}^r\,\ep_i^{m_i}$, set
$$f_{(\de)}:=f|_{X_{(\de)}}:X_{(\de)}\to\De_{\de}\q\h{with}\q X_{(\de)}:=X\caps f^{-1}(\De_{\de}).$$
This is better than the Milnor fibration, since we can apply the K\"unneth formula to the fibers. The canonical morphism (1.5.5) is constructed using this fibration together with the canonical quasi-isomorphism
$$f_U^{-1}\OO_{\De^*}\sot_{\C}L\1[d_X{-}1]\simto\DR_{X/\De}^{\log}(\Lc_X)|_U.$$
\sk
Using a filtration on $L$, the assertion is then reduced to the case
$${\rm rank}\,L=1,$$
since the local monodromies commute with each other in the case of the complement of a divisor with normal crossings. (Recall that $(\De^*)^r$ has a deformation retract to a real torus and its fundamental group is an abelian group $\Z^r$. This real torus is the product of the topological cycles $\ga_j$ in the proof of Prop.\,1.3d.)
\sk
Let $\la_i$ be the eigenvalue of the local monodromy of $L$ around $\{z_i\eq 1\}$ ($i\ins[1,n]$). Let $\al_i\ins[0,1)$ with $\la\eq e^{-2\pi\sqrt{-1}\1\al_i}$.
By definition $\DR_{X/\De}^{\log}(\Lc_X)_x$ can be identified with the Koszul complex for the actions of $\eta_i$ on $\OO_{\!X,x}\eq\C\{z_1,\dots,z_{d_X}\}$ ($i\ins[2,d_X]$) with
$$\eta_i:=\begin{cases}(z_i\dd_{z_i}{+}\,\al_i)-\tfrac{m_i}{m_1}(z_1\dd_{z_1}{+}\,\al_1)&(i\ins[2,r])\\ z_i\dd_{z_i}{+}\,\al_i&(i\ins[r{+}1,n])\\ \dd_{z_i}&(i\ins[n{+}1,d_X]),\end{cases}$$
see \cite[Prop.\,1.13]{St} for the case $L\eq\C_U$. Recall that
$$\aligned&\Om^p_{X/\De}\lYr=\Om^p_X\lYr/(\ddd f\1\!/\1\!f){\wedge}\Om^{p-1}_X\lYr,\\&\h{with}\q\q\ddd f\1\!/\1\!f=\msum_{i=1}^r\,m_i\1\ddd z_i/z_i,\endaligned$$
and $\Om\ub_X\lYr_x$ (resp. $\Om\ub_{X/\De}\lYr_x$) has a free basis consisting of exterior products of $\ddd z_i/z_i$ for $i\ins[1,n]$ (resp. $i\ins[2,n]$) and $\ddd z_i$ for $i\ins[n{+}1,d_X]$.
\sk
We then see that the Koszul complex is acyclic if $\al_i/m_i\ne\al_j/m_j$ for some $i,j\ins[1,r]$ or $\al_i\ne 1$ for some $i\ins[r{+}1,n]$. The corresponding vanishing holds for the right-hand side of (1.5.5) using the K\"unneth formula. Here we may assume that $c\eq1$ and $\ep_i\eq{+}\infty$ ($i\in[1,d_X]$), that is, $X\eq\C^{d_X}$ for the calculation of the right-hand side. We use a partial compactification $\C^{d_X}\subset\PP^1{\times}\C^{d_X-1}$, and consider the projection to the divisor at infinity, where $z_1^{-1},z_2,\dots,z_{d_X}$ give local coordinates around the divisor.
\sk
We may thus assume
$$\al_i/m_i\eq\al_j/m_j\q(i,j\in[1,r]),\q\al_i\eq0\q(i\in[r{+}1,n]).
\leqno(1.5.6)$$
Here we may also assume $n\eq d_X$, and then $r\eq n$ using the triple complex structure of the Koszul complex together with the K\"unneth formula for $U_c\eq U'_c{\times}U''$, where $U_c\eq U\caps X_c$ and $U''\eq(\De^*)^{n-r}$, see also the proof of Prop.\,1.3d. We thus get the shifted Koszul complex for the action of
$$z_i\dd_{z_i}{-}\,(m_i/m_1)z_1\dd_{z_1}\q(i\ins[2,d_X]),$$
using the bases consisting of exterior products of $\ddd z_i/z_i$ ($i\ins[2,d_X]$).
\sk
This Koszul complex is quasi-isomorphic to the Koszul subcomplex for the {\it zero actions\1} on $\C\{t'\}\sst\C\{z_1,\dots,z_{d_X}\}$, where $t':=f':=\mprod_{i=1}^{d_X}\,z_i^{m'_i}$ with $m'_i:=m_i/m$ and $m:={\rm GCD}(m_i)$, see also the proof of \cite[Prop.\,1.13]{St}. Note that $\C\{t'\}$ is a free $\C\{t\}$-module of rank $m$. Here $t$ is identified with
$$f^*t=f=\mprod_{i=1}^{d_X}\,z_i^{m_i}=f'{}^m=t'{}^m.$$
\sk
We can verify that the action of $t\dd_t$ on an exterior product of $\ddd z_i/z_i$ for $i\ins[2,d_X]$ (which is a member of the basis of the relative logarithmic de Rham complex) multiplied by $t'{}^a$ ($a\ins[0,m{-}1]$) is given by
$$(\al_1{+}\,a\1m'_1)/m_1=\al_1/m_1\pl a/m\,\in\,[0,1),
\leqno(1.5.7)$$
using the short exact sequence (1.5.5) (whose connection morphism is used for the definition of the action of $t\dd_t$). Here we look at the coefficients of the exterior product of the above exterior product of $\ddd z_i/z_i$ ($i\ins[2,d_X]$) with $\ddd z_1/z_1$.
\sk
The above calculation is compatible with the corresponding topological calculation for the right-hand side of (1.5.5). Note that the hypothesis $\al_i/m_i\eq\al_j/m_j$ in (1.5.6) means that the local system $L$ is the pullback of a local system on $\De^*$ by $f_U$. So its restriction on $U_c$ is constant. Assuming $U\eq(\C^*)^{d_X}$ as above, the number of connected components of $U_c$ is $m$, and each connected component is isomorphic to an affine torus, since it is an unramified covering of $(\C^*)^{d_X-1}$. So both sides of (1.5.5) have the same dimension.
\sk
To show that (1.5.5) is an isomorphism, it is then sufficient to show that the integration of the pullback of an exterior product of $\ddd z_i/z_i$ ($i\ins[2,d_X]$) over some topological cycle of each connected component of $U_c$ does not vanish. Here we can consider the integration of the form over the image of the topological cycle in $(\C^*)^{d_X-1}$. However, any finite unramified covering of affine tori induces an isomorphism of cohomology groups with $\Q$-coefficients. (Indeed, it induces an inclusion of their fundamental groups, which are isomorphic to their first homology groups, with a finite cokernel.) So the assertion follows. This finishes the proof of Prop.\,1.5a.
\sk
The following gives a generalization of a variant of \cite[Prop.\,2.4]{gau} (where the eigenvalues of the residue of $t\dd_t$ are contained in $(-1,0]$, since the logarithmic forms with poles along $X_0$ is not used there).
\msn
{\bf Corollary\hs1.5a.} {\it Let $f:X\tos\De$ be as in the beginning of {\bf 1.5}. Assume $f$ is proper. Put
$$f_*^{\1\DR^{\log}}\Lc_X:=\RR f_*\DR_X^{\log}\bl(\Lc_X[\dd_tt]\1\de(f{-}t)\br).$$
Then the $\Hc^jf_*^{\1\DR^{\log}}\Lc_X$ are free of finite rank over $\OO_{\De}$, and the real part of any eigenvalue of the residue of $\dd_tt$ is contained in $[0,1)$, shrinking $\De$ if necessary.}
\msn
{\it Proof.} We may assume that $Y\stm X_0$ is a divisor with relative normal crossings over $\De^*$ shrinking $\De$ if necessary.
As is shown in the proof of Prop.\,1.3e, the restrictions of $\Hc^jf_*^{\1\DR^{\log}}\Lc_X$ to $\De^*$ are locally free of finite rank over $\OO_{\De^*}$. Moreover the $\Hc^jf_*^{\1\DR^{\log}}\Lc_X$ are coherent over $\OO_{\De}$ by the Grauert coherence theorem, see for instance \cite{Dou}. (Note that the quasi-isomorphisms in (1.5.2--3) hold as complexes of $f^{-1}\OO_{\De}$-modules.)
\sk
Since $f$ is proper, we have the isomorphisms
$$\Hc^j\bl(f_*^{\1\DR^{\log}}\Lc_X\br){}_0=H^j\bl(X_0,\DR_X^{\log}\bl(\Lc_X[\dd_tt]\1\de(f{-}t)\br)|_{X_0}\br),$$
together with a Leray-type spectral sequence
$$E_2^{p,q}=H^p\bl(X_0,\Hc^q\1\DR_X^{\log}\bl(\Lc_X[\dd_tt]\1\de(f{-}t)\br)|_{X_0}\br)\Longrightarrow\Hc^{p+q}\bl(f_*^{\1\DR^{\log}}\Lc_X\br){}_0,$$
which is defined in the category of $\C\{t\}$-modules with the action of $\dd_tt$. In order to show the assertions of Cor.\,1.5a, it is then enough to prove that the $E_2$-terms are finite free $\C\{t\}$-modules and the real part of any eigenvalue of the residue of $\dd_tt$ is contained in $[0,1)$, see also Rem.\,1.5b. However, these follow from Prop.\,1.5a. So the assertions are proved. This finishes the proof of Cor.\,1.5a.
\msn
{\bf Corollary\hs1.5b.} {\it In the notation and assumption of {\rm Cor.\,1.5a}, the $\D_{\De}$-modules $\Hc^j\!f_*\uD\!\M_X(L)$ are regular holonomic of normal crossing type shrinking $\De$ if necessary.}
\msn
{\it Proof.} By the construction at the beginning of {\bf 1.5}, we have a natural inclusion
$$\DR_X^{\log}\bl(\Lc_X[\dd_tt]\1\de(f{-}t)\br)\into\DR_X^{\log}\bl(\Lc_X(*Y)[\dd_tt]\1\de(f{-}t)\br),$$
and the target is isomorphic to $\DR_X\bl(\M_X(L)[\dd_t]\1\de(f{-}t)\br)$. So we get the canonical morphisms
$$\Hc^jf_*^{\1\DR^{\log}}\Lc_X\to\Hc^jf_*\uD\M_X(L)\q(j\ins\Z),$$
inducing the isomorphisms
$$\Hc^jf_*^{\1\DR^{\log}}\Lc_X[t^{-1}]\simto\Hc^jf_*\uD\M_X(L),$$
where $f_*^{\1\DR^{\log}}\Lc_X$ is as in Cor.\,1.5a. The assertion then follows from Cor.\,1.5a. This finishes the proof of Cor.\,1.5b.
\msn
{\bf Corollary\hs1.5c.} {\it Let $f:X\to S$ be a proper morphism of complex manifolds. Let $Y\sst X$, $D\sst S$ be divisors with normal crossings such that $X_D\,{:=}\,f^{-1}(D)\sst Y$, $f$ is smooth over $S\stm D$, and $Y\stm X_D$ is a divisor with relative normal crossings over $S\stm D$. Let $L$ be a $\C$-local system on $U\,{:=}\,X\stm Y$ with $\M_X(L)$ as in {\rm Prop.\,1.3b}. Then the $\Hc^jf_*\uD\M_X(L)$ are meromorphic Deligne extensions of local systems on $S\stm D$.}
\msn
{\it Proof.} By Prop.\,1.3e, the restrictions of $\Hc^jf_*\uD\M_X(L)$ over $S\stm D$ are locally free of finite rank over $\OO_{S\setminus D}$. By {\bf A.6} below, the $\Hc^jf_*\uD\M_X(L)$ are coherent $\OO_S$-modules, since the hypothesis on the existence of a coherent $\OX$-submodule generating the $\D_X$-module over $\D_X$ is satisfied by using the Deligne extension, see Rem.\,1.3e. By Prop.\,1.3, the assertion is then reduced to the $D$ {\it smooth\1} case (shrinking $S$). This case may be viewed as the situation of Cor.\,1.5b with parameters.
\sk
We proceed by induction on $d_S$. The assertion follows from Cor.\,1.5b if $d_S\eq 1$. Assume $d_S\gess 2$. Let $t_1,\dots,t_{d_S}$ be local coordinates of $S$ such that $D\eq\{t_1\eq0\}$ locally. Set $f_i\,{:=}\,f^*t_i$ ($i\ins[1,d_S]$). We may assume that $(f_2,\dots,f_{d_S}):X\to\C^{d_S}$ is a smooth morphism shrinking $S$ if necessary (using a Bertini-type theorem).
\sk
Generalizing the construction at the beginning of {\bf 1.5}, consider the complex
$$\Cc_f\ub:=\DR_X^{\log}\bl(\Lc_X[\dd_{t_1}t_1,\dd_{t_2},\dots,\dd_{t_{d_S}}]\1\de(f_1{-}t_1)\cdots\de(f_{d_S}{-}t_{d_S})\br),$$
which is a subcomplex of
$$\DR_X\bl(\M_X(L)[\dd_{t_1},\dd_{t_2},\dots,\dd_{t_{d_S}}]\1\de(f_1{-}t_1)\cdots\de(f_{d_S}{-}t_{d_S})\br).$$
By an argument similar to the proof of Prop.\,1.3e, we can verify that the former complex is {\it quasi-isomorphic\1} as a complex of $f^{-1}\OO_S$-modules to the subcomplex
$$\Cc_f^{\prime\ssb}:=\DR_{X,\{d_S\}}^{\log}\bl(\Lc_X[\dd_{t_1}t_1,\dd_{t_2},\dots,\dd_{t_{d_S-1}}]\1\de(f_1{-}t_1)\cdots\de(f_{d_S-1}{-}t_{d_S-1})\br),$$
where $\DR_{X,\{d_S\}}^{\log}$ means that we replace the $\Om_X^p\lYr$ with ${\rm Ker}\bl(\ddd f_{d_S}\wedge:\Om_X^p\lYr\to\Om_X^{p+1}\lYr\br)$. Note that this complex is further quasi-isomorphic as a complex of $f^{-1}(\OO_S)$-modules to
$$\DR^{\log}_{X/S}(\Lc_X),$$
and its direct image sheaves are {\it coherent\1} by the Grauert theorem, although the relation to the Gauss-Manin connection on these does not seem trivial.
\sk
By inductive hypothesis, we may assume that the assertion holds for the restriction
$$f_{(c)}\,{:=}\,f|_{X_{(c)}}:X_{(c)}\tos S_{(c)}\q\h{with}\q X_{(c)}\,{:=}\,f^{-1}(S_{(c)}),$$
and $S_{(c)}\sst S$ is defined by $t_{d_S}\eq c$ for $c\ins\De_{\de}$. There is a short exact sequence of complexes
$$0\to\Cc_f^{\prime\ssb}\buildrel{\!\mu_c}\over\longrightarrow\Cc_f^{\prime\ssb}\to\Cc_{f_{(c)}}\ub\to 0,$$
where the morphism $\mu_c$ is defined by multiplication by $f_{d_S}{-}\1c$. This gives the long exact sequence of coherent sheaves for $c\ins\De_{\de}$\,:
$$\to R^j\!f_*\Cc_f^{\prime\ssb}\buildrel{\!\mu_c^j}\over\longrightarrow R^j\!f_*\Cc_f^{\prime\ssb}\to R^j(f_{(c)})_*\Cc_{f_{(c)}}\ub\to R^{j+1}f_*\Cc_f^{\prime\ssb}\buildrel{\mu_c^{j+1}}\over\longrightarrow.$$
From this, we can deduce by induction on $d_S$ that the $R^j\!f_*\Cc_f\ub$ are locally free sheaves of finite rank, since the injectivity of $\mu_c^j$ on the direct image sheaves $R^j\!f_*\Cc_f^{\prime\ssb}$ follows from the freeness of the $R^j(f_{(c)})_*\Cc_{f_{(c)}}\ub$ by {\it decreasing\1} induction on $j$ using Nakayama's lemma. Indeed, setting
$$\rho_j:={\rm rank}\,\Hc^jf_*\uD\M_X(L)|_{S\setminus D},$$
the sequence implies that the coherent sheaf $R^j\!f_*\Cc_f^{\prime\ssb}$ is locally generated by $\rho_j$ element over $\OO_S$ (and its restriction to $S\stm D$ is locally free of rank $\rho_j$). Consider the kernel of the surjection
$$\buildrel{\rho_j}\over{\mopl}\OO_S\to R^j\!f_*\Cc_f^{\prime\ssb}.$$
This cannot be supported on $D$. So the freeness of $R^j\!f_*\Cc_f^{\prime\ssb}$ follows.
\sk
We can verify that the real part of any eigenvalue of the residue of $\dd_{t_1}t_1$ is contained in $[0,1)$ by restricting to a curve defined by $t_i\eq c_i$ ($i\ins[2,d_S]$). (Here we can also use {\it strict Nilsson class functions.}) We employ modules over the differential operators $\D_{X/T^k}$, $\D_{S/T^k}$ which are subrings of $\D_X$, $\D_S$, where $T^k:=\C^k$ with coordinates $t_i $ ($i\ins[d_S{-}k{+}1,d_S]$) for $k\ins[1,d_S{-}1]$. These modules are equivalent to $\OO$-modules endowed with an integrable {\it relative\1} connection, where relative differential forms $\Om^{\ssb}_{X/T^k}$, $\Om^{\ssb}_{S/T^k}$ are used. We can show that the passage to relative $\D$-modules (forgetting the action of certain differential operators) commutes with the direct image by $f$, using the factorization $X\tos X{\times}S\tos S$, and applying an argument similar to the proof of Prop.\,1.3e so that the relative differential forms $\Om^{\ssb}_{X\times S/S}$ in the definition of the direct image of $\D$-modules under the projection to $S$ are replaced by $\Om^{\ssb}_{X\times_{T^k}S/S}$ by increasing induction on $k$. (Note that $X{\times}_{T^k}S\to S$ is the base change of $X\tos T^k$ by $S\tos T^k$.) The assertion for $d_S\eq1$ is proved in Cor.\,1.5a. This finishes the proof of Cor.\,1.5c.
\bs\bs
\vbox{\centerline{\bf 2. Formalism of regular holonomic $\D$-modules}
\bsn
In this section we give a definition of regular holonomic $\D$-modules using algebraic local cohomology and meromorphic Deligne extensions in {\bf 2.1}, and prove some of their properties in {\bf 2.2}.}
\msn
{\bf 2.1.~Regular holonomic $\D$-modules.} Let $X$ be a smooth complex variety. We say that a holonomic left $\D_X$-module $\M$ is {\it regular holonomic\1} at $x\ins X$ if the following conditions are satisfied by induction on $d_Z$ with $Z:={\rm Supp}\,\M$ {\it replacing $X$ with a sufficiently small open neighborhood of $x$ in $X$}.
\sk
(i) There is a closed subvariety $(Y,x)\sst (X,x)$ such that $Z\stm Y$ is purely $d_Z$-dimensional and smooth with $d_{Y\cap Z}\,{<}\,d_Z$, the restriction of $\DRt_X(\M)[-d_Z]$ to $Z\stm Y$ is a local system, and moreover $\Hc^0_{[Y]}\M$ is {\it regular holonomic.}
\sk
(ii) For each $d_Z$-dimensional irreducible component $Z_i$ of $Z$, there is a desingularization $\pi'_i:\Zt_i\tos Z_i$ such that $\Yt_i:=\pi_i^{\prime-1}(Y_i)$ with $Y_i:=Y\caps Z_i$ is a divisor with normal crossings on $\Zt_i$, $\pi'_i$ induces an isomorphism over $Z_i\stm Y_i$, and there is an isomorphism of $\D_X$-modules
$$\Hc^0_{[X|Y]}\M\cong\mopl_i\,\Hc^0(\pi_i)_*\uD\Mt_i.
\leqno(2.1.1)$$
Here $\pi_i$ is the composition of $\pi'_i$ with the inclusion $Z_i\into X$, and $\Mt_i$ is the {\it meromorphic Deligne extension\1} on $\Zt_i$ of the local system $\DRt_X(\M)[-d_Z]|_{Z_i\setminus Y_i}$ viewed as a $\D_{\Zt_i}$-module, see Prop.\,1.3b. (Note that $\Hc^0(\pi_i)_*\uD$ is the {\it composition\1} of $(\pi_i)_*\uD$ and the cohomology sheaf functor $\Hc^0$.)
\sk
Holonomic $\D$-modules with discrete supports are regular holonomic, since any holonomic $\D_X$-module with $X\eq{\rm pt}$ is regular holonomic of normal crossing type where $Y\eq\emptyset$.
\sk
We will denote by $\Msf_{\rh}(\D_X)\sst \Msf(\D_X)$ the full subcategory consisting of regular holonomic left $\D_X$-modules. The derived category of bounded complexes of left $\D_X$-module with regular holonomic cohomology sheaves will be denote by $\Dsf^b_{\rh}(\D_X)$.
\msn
{\bf Remark\hs2.1a.} Under condition~(i) (with last condition on regular holonomicity forgotten), there is a canonical direct sum decomposition
$$\Hc^0_{[X|Y]}\M\cong\mopl_i\,\M_i,
\leqno(2.1.2)$$
using Rem.\,1.2a, where the $\M_i$ are $\D_X$-modules with ${\rm Supp}\,\M_i\eq Z_i$. Indeed, one can put
$$\M_i:=\Hc_{[X|Y]}^0\Hc_{[Z_i]}^0\M.$$
Applying $\Hc_{[X|Y]}^0\Hc_{[Z_i]}^0$ to (2.1.1), we also get that
$$\M_i=\Hc^0(\pi_i)_*\uD\Mt_i.$$
So the $\M_i$ are regular holonomic by the above definition (since $\Hc^0_{[Y]}\M_i\eq0$, $\Hc^0_{[X|Y]}\M_i\eq\M_i$). The $\M_i$ are then called the {\it meromorphic regular holonomic extensions\1} of the local systems $L_i:=\DRt_X(\M)[-d_Z]|_{Z_i\setminus Y_i}$. The smallest regular holonomic $\D_X$-submodule of $\M_i$ whose quotient is supported in $Z_i\stm Y_i$ will be called the {\it minimal regular holonomic extensions\1} of the local systems $L_i$, see also Rem.\,1.2d.
\msn
{\bf Remark\hs2.1b.} Condition~(2.1.1) is equivalent to that $\Hc^0_{[\Zt_i|\Yt_i]}(\pi_i)_{\D}^!\Hc^0_{[X|Y]}\M$ is isomorphic to the meromorphic Deligne extension $\Mt_i$. This can be verified using the canonical morphisms
$$(\pi_i)_*^{\D}\ssc(\pi_i)^!_{\D}\to{\rm id},\q{\rm id}\to(\pi_i)^!_{\D}\ssc(\pi_i)_*^{\D},$$
which follow from the adjunction relation between $(\pi_i)_*^{\D}$ and $(\pi_i)^!_{\D}$.
\sk
It does not seem, however, quite easy to prove this adjunction relation. Alternatively, one can use a composition of functors in Prop.\,1.4a--b instead of $\Hc^0_{[\Zt_i|\Yt_i]}(\pi_i)_{\D}^!$.
\msn
{\bf Proposition\hs2.1.} {\it Condition~{\rm(ii)} is independent of the choice of $\pi_i:\Zt_i\tos Z_i$. Conditions~{\rm(i)} and {\rm(ii)} are independent of the choice of $Y$ satisfying the conditions in {\rm(i)} except the last condition on regular holonomicity of $\Hc^0_{[Y]}\M$.}
\msn
{\it Proof.} For the first assertion, we apply Prop.\,1.2 and 1.3b to a morphism between two desingularizations, which induces an isomorphism over an open subset whose complement has codimension at least 2. Then the assertion follows from Prop.\,1.3c.
\sk
The second assertion for condition~(i) is clear. As to condition~(ii), assume there are two subspaces $Y_1,Y_2$ satisfying the conditions of (i) except the last one. Here we may assume $Y_1\sst Y_2$ (taking the intersection). For $a\eq1,2$, there are desingularizations $\pi_{i,a}:\Zt_{i,a}\tos Z_i$ together with the divisors with normal crossings
$$\Yt_{i,a}:=\pi_{i,a}^{-1}(Y_{i,a})\sst\Zt_{1,a}.$$
Here $Y_{i,a}:= Z_i\caps Y_a$, and $\pi_{i,a}$ induces an isomorphism over $Z_i\stm Y_{i,a}$. By the first assertion we may assume that $\pi_{i,2}=\pi_{i,1}\ssc\pi_{i,0}$ using an embedded resolution
$$\pi_{i,0}:\bl(\Zt_{i,2},\Yt_{i,2}\br)\to\bl(\Zt_{i,1},\Yt^+_{i,1}\br)\q\h{with}\q\Yt^+_{i,1}:=\pi_{i,1}^{-1}(Y_{i,2}).$$
\sk
Set
$$\Yt'_{i,2}:=\pi_{i,0}^{-1}(\Yt_{i,1})\sst\Zt_{i,2}.$$
This is a divisor with normal crossings, since it is a divisor contained in $\Yt_{i,2}\sst\Zt_{i,2}$. (Here we {\it cannot\1} assume $\Zt_{i,1}\eq\Zt_{i,2}$ replacing $\Zt_{i,1}$ with $\Zt_{i,2}$ and using the independence of $\Zt_i$ proved above, since $\pi_{i,2}$ does not necessarily induce an isomorphism over $Z_i\stm Y_{i,1}$.)
\sk
Let $\Mt'_{i,2}$ be the $\D_{\Zt_{i,2}}$-module which is the meromorphic Deligne extension of the pullback of the local system $\DRt_X(\M)[-d_Z]|_{Z_i\setminus Y_{i,1}}$. We have also the $\D_{\Zt_{i,a}}$-module $\Mt_{i,a}$ which is the meromorphic Deligne extension of the pullback of the local system $\DRt_X(\M)[-d_Z]|_{Z_i\setminus Y_{i,a}}$ ($a\eq1,2$) as in condition~(ii).
\sk
By the uniqueness of Deligne extension, we have the natural isomorphism
$$\RR\Ga_{[\Zt_{i,2}|\Yt_{i,2}]}\Mt'_{i,2}\eq\Mt_{i,2}.$$
Applying an argument similar to the proof of Prop.\,1.3c (using Hartogs extension theorem), there is a canonical morphism
$$\Hc^0(\pi_{i,0})_*\uD\Mt'_{i,2}\to\Mt_{i,1},$$
inducing an isomorphism on the complement of $\Yt^+_{i,1}$ (since $\pi_{i,0}$ induces an isomorphism over the complement of a closed subvariety with codimension at least 2). The cohomology sheaves of the mapping cone of this morphism is supported in $\Yt^+_{i,1}$, hence it induces an isomorphism
$$(\pi_{i,0})_*\uD\Mt_{i,2}=\RR\Ga_{[\Zt_{i,1}|\Yt^+_{i,1}]}\Mt_{i,1},$$
by Rem.\,1.2a and Prop.\,1.2 (commutativity of local cohomology and direct image). By the latter proposition (applied to $(\pi_{i,1})_*\uD$), we then get the morphism
$$\iota_i:(\pi_{i,1})_*\uD\Mt_{i,1}\to\RR\Ga_{[X|Y_2]}(\pi_{i,1})_*\uD\Mt_{i,1}=(\pi_{i,2})_*\uD\Mt_{i,2},$$
inducing an isomorphism on the complement of $Y$.
\sk
These imply the desired independence. Indeed, to show that the isomorphism (2.1.1) for $Y_2$ implies the one for $Y_1$ (where the converse is easy), we consider the kernel of the composition
$$\aligned&\Hc^0_{[X|Y_1]}\M\to\Hc^0_{[X|Y_2]}\M\cong\mopl_i\,\Hc^0(\pi_{i,2})_*\uD\Mt_{i,2}\\ &\to\mopl_i\,{\rm Coker}\bl(\Hc^0\iota_i:\Hc^0(\pi_{i,1})_*\uD\Mt_{i,1}\into\Hc^0(\pi_{i,2})_*\uD\Mt_{i,2}\br),\endaligned$$
which has a morphism to $\mopl_i\,\Hc^0(\pi_{i,1})_*\uD\Mt_{i,1}$. By Rem.\,1.2a we then get the isomorphism
$$\Hc^0_{[X|Y_1]}\M=\mopl_i\,\Hc^0(\pi_{i,1})_*\uD\Mt_{i,1},$$
since the support of the image of the above composition is contained in $Y_1$. (Note that the restriction of $\DRt_X(\M)$ to the complement of $Y_1$ is a shifted local system.) This finishes the proof of Prop.\,2.1.
\sk
Prop.\,2.1 implies the following.
\msn
{\bf Corollary\hs2.1a.} (1) {\it The subcategory $\Msf_{\rh}(\D_X)$ is closed by taking subobjects and extensions $($in particular, direct sum$)$ in $\Msf_{\rm coh}(\D_X)$, that is, for a short exact sequence
$$0\tos\M'\tos\M\tos\M''\tos0\q\h{in}\,\,\,\,\,\Msf_{\rm coh}(\D_X),$$
we have $\M'\ins \Msf_{\rh}(\D_X)$ if $\M\ins \Msf_{\rh}(\D_X)$, and $\M\ins \Msf_{\rh}(\D_X)$ if $\M',\M''\ins \Msf_{\rh}(\D_X)$.
\sk
{\rm (2)} The subcategory $\Msf_{\rh}(\D_X)$ is stable by $\Hc^0_{[X|D]}$ for any divisor $D\sst X$.}
\ms
(The stability by taking quotient objects in $\Msf_{\rm coh}(\D_X)$ will be proved in Prop.\,2.2 below.)
\msn
{\it Proof.} We prove the assertions by increasing induction on $d_Z$ with $Z:={\rm Supp}\,\M$. It is clear for $d_Z\eq0$. Assume there is a short exact sequence $0\tos\M'\tos\M\tos\M''\tos 0$ as above.
\sk
Assume first that $\M',\M''$ are regular holonomic. Set $Z':={\rm Supp}\,\M'$, $Z'':={\rm Supp}\,\M''$. By Prop.\,2.1 we may assume that condition~(i) for $\M',\M''$ is satisfied for the same {\it divisor\1} $Y\subset X$ shrinking $X$ if necessary. We may also assume that the desingularization of an irreducible component $Z'_i\subset Z'$ and that of $Z''_j\subset Z''$ are the same in the case $Z'_i\eq Z''_j$. Condition~(i) is then satisfied for $\M$ by inductive hypothesis together with the {\it left exactness\1} of $\Hc_{[Y]}^0$. For condition~(ii) we apply Prop.\,1.2 and 1.3a--b together with the exactness of $\Hc_{[X|Y]}^0$ and also Rem.\,2.1b. So $\M$ is regular holonomic. The argument is similar and much easier in the case $\M$ is assumed to be regular holonomic.
\sk
For the proof of the stability by $\Hc^0_{[X|D]}$ in (2), we have the following exact sequence using the {\it left exactness\1} of the functor $\Hc^0_{[X|D]}$\,:
$$0\to\Hc^0_{[X|D]}\Hc_{[Y]}^0\M\to\Hc^0_{[X|D]}\to\Hc_{[X|D]}^0\Hc_{[X|Y]}^0\M.$$
The assertion (2) then follows from the inductive hypothesis and the stability by subobjects and extensions proved above together with the isomorphism (2.1.1) and Prop.\,1.2, Rem.\,1.2c. Here we replace $\Zt_i$ with the embedded resolution of $Y'_i:=(Y\cup D)\cap Z_i\subset Z_i$, and irreducible components $Z_i$ contained in $D$ are neglected. This finishes the proof of Cor.\,2.1a.
\sk
From Cor.\,2.1a, we can deduce the following.
\msn
{\bf Corollary\hs2.1b.} {\it Any regular holonomic $\D$-module of normal crossing type on a complex manifold $($see {\bf 1.3}$)$ is regular holonomic.}
\msn
{\it Proof.} Since the assertion is local, we may assume $X$ is a polydisk $\De^{d_X}$ with local coordinates $z\eq(z_1,\dots,z_{d_X})$. By Cor.\,2.1a, we may assume that $\M_x\ins \Msf_{\rhnc}^{(z)}$ is {\it simple.} Then we have $\M_{x,z}^{(\nu)}\cong\C$ for some $\nu\ins\Xi^{d_X}$, and $\M_{x,z}^{(\nu')}\eq 0$ for any $\nu'\ins\Xi^{d_X}{\setminus}\,\{\nu\}$ as in the proof of Prop.\,1.3a, where we may assume moreover $\nu_i\ne{-1}$ for any $i$. Then the restriction of $\DRt_X(\M)$ to $X\stm Y=(\De^*)^{d_X}$ is a shifted local system of rank 1. The assertion then follows from Rem.\,1.3c.
\ms
Prop.\,2.1 implies also the following.
\msn
{\bf Corollary\hs2.1c.} {\it Let $\M$ be a holonomic $\D_X$-module. Let $\{Z_k\}_{k\in\N}$ be an increasing sequence of closed subvarieties associated with a Whitney stratification of ${\rm Supp}\,\M$, that is, setting $S_k:=Z_k\stm Z_{k-1}$ $(\forall\,k\ins\N)$ with $Z_{-1}:=\emptyset$, $\{S_k\}_{k\in\N}$ is a Whitney stratification of ${\rm Supp}\,\M$ with $\dim S_k\eq k$. Assume the $\Hc^j\DRt_X(\M)|_{S_k}$ are local systems $(j\ins\Z,\,k\ins\N)$.  Then $\M$ is regular holonomic if and only if the following condition holds locally on $X:$
\sk
The $\Hc^0_{[X|Y_k]}\Hc^0_{[Z_k]}\M$ correspond to regular meromorphic connections via {\rm Prop.\,1.4b}, where the $Y_k\sst X$ are locally defined divisors containing $Z_{k-1}$ with $\dim Y_k\caps Z_k<k$ $(\forall\,k\,{>}\,0)$.}
\msn
{\bf Corollary\hs2.1d.} {\it For a simple regular holonomic $\D_X$-module $\M$, set $Z:={\rm Supp}\,\M$, and let $Y\sst Z$ be the largest open subset such that $\DRt_X(\M)[-d_Z]|_{Z\setminus Y}$ is a local system, denoted by $L$. Then $Z$ is globally irreducible, $Z\stm Y$ is connected, $L$ is simple, and $\M$ is the minimal regular holonomic extension of $L$, see {\rm Rem.\,1.2d} and {\rm 2.1a} for minimal extensions.}
\msn
{\it Proof.} Let $Z_i$ be the global irreducible components of $Z$. The argument in the beginning of {\bf 2.1} can be applied globally, and we get the inclusion
$$\M\into\mopl_i\,\M_i,$$
with $\M_i$ as in Rem.\,2.1a, since $\Hc^0_{[Y]}\M\eq0$ by the simplicity of $\M$. Choosing a total order of the index set of $\M_i$, we get a finite filtration on the direct sum whose graded quotients are $\M_i$. Restricting to the complement of $Y$, this implies that $Z$ must be globally irreducible, that is, $Z\stm Y$ is connected, and we get the inclusion $\M\into\M_1$. If $L$ is not simple, we have a nontrivial filtration of $\M_1$ whose restriction to the complement of $Y$ is also nontrivial. Considering its induced filtration on $\M$, this gives a contradiction, since $\M|_{X\stm Y}\eq\M_1|_{X\stm Y}$. It is then easy to see that $\M$ is the minimal regular holonomic extension of $L$, see Rem.\,1.2d and 2.1a. This finishes the proof of Cor.\,2.1d.
\msn
{\bf 2.2.~Proof of Theorem~1.} We first show the following.
\msn
{\bf Proposition\hs2.2.} (1) {\it For $\M\ins \Msf_{\rh}(\D_X)$ and a surjection $\M\tos\M''$ with $\M''\ins \Msf_{\rm coh}(\D_X)$, we have $\M''\ins \Msf_{\rh}(\D_X)$.
\sk
{\rm (2)} For a proper morphism of smooth complex varieties $f\,{:}\,X\tos S$ and $\M\ins \Msf_{\rh}(\D_X)$, we have $\Hc^j\!f\uD_*\M\ins \Msf_{\rh}(\D_S)$ $(j\ins\Z)$.}
\msn
{\it Proof.} We prove these assertions by induction on $d_{\!\M}\,{:=}\,\dim{\rm Supp}\,\M$. The assertions are clear if $d_{\!\M}\eq 0$. Let $d$ be a positive integer. Assume the assertions are proved for $d_{\!\M}<d$. We first prove (1) for $d_{\!\M}\eq d$. Assume there is a short exact sequence of $\D_X$-modules
$$0\to\M'\to\M\to\M''\to0,$$
with $\M',\M$ regular holonomic (using Cor.\,2.1a). Let $Y\sst X$ be a {\it divisor\1} such that conditions~(i) and (ii) are satisfied for $\M$ (shrinking $X$ if necessary). There are exact sequences
$$\aligned&0\to\Hc_{[Y]}^0\M'\to\Hc_{[Y]}^0\M\to\Hc_{[Y]}^0\M''\to\Hc_{[Y]}^1\M',\\&0\to\Hc_{[X|Y]}^0\M'\to\Hc_{[X|Y]}^0\M\to\Hc_{[X|Y]}^0\M''\to0.\endaligned$$
Condition~(ii) for $\M''$ is then satisfied using Prop.\,1.3a and 1.4a--b.
\sk
For the proof of condition~(i), it is enough to show that $\Hc_{[Y]}^1\M'$ is regular holonomic using Cor.\,2.1a and the inductive hypothesis. Here we may replace $\M'$ with $\M$ to simplify the notation (since $\M$ is any regular holonomic $\D_X$-module). By the exact sequence
$$0\to\Hc_{[Y]}^0\M\to\M\to\Hc_{[X|Y]}^0\M\to\Hc_{[Y]}^1\M\to 0,$$
it is sufficient to show that any coherent quotient $\D_X$-module of $\Hc_{[X|Y]}^0\M$ whose support is contained in $Y$ is regular holonomic.
\sk
Let $\pi_i:\Zt_i\to Z_i$ and $\Mt$ be as in condition~(ii). Let $\Mt'_i\sst\Mt_i$ be the smallest $\D_{\Zt_i}$-submodule such that ${\rm Supp}\,(\Mt_i/\Mt'_i)\sst\Yt_i$. This is the image of the canonical morphism
$$^{\vee}\!\!\Hc^0_{[\Zt_i|\Yt_i]}\Mt^i:=\DD\bl(\Hc^0_{[\Zt_i|\Yt_i]}\DD\Mt_i\br)\to\Mt_i,$$
see Rem.\,1.2d. Here $\DD$ denotes the dual functor.
\sk
Set $\Mt_i'':=\Mt_i/\Mt'_i$ so that we have the exact sequence
$$0\to\Hc^{-1}(\pi_i)_*^{\D}\Mt''_i\to\Hc^0(\pi_i)_*^{\D}\Mt'_i\to\Hc^0(\pi_i)_*^{\D}\Mt_i\to\Hc^0(\pi_i)_*^{\D}\Mt''_i.$$
Here $\Hc^{-1}(\pi_i)_*^{\D}\Mt=0$ by Prop.\,1.2 and Rem.\,1.2b. Note that $\Mt_i\eq\Hc_{[\Zt_i|\Yt_i]}\Mt'_i$.
\sk
Using (2.1.1) and the inductive hypothesis, it is then enough to show that any coherent quotient $\D_X$-module of $\Hc^0(\pi_i)_*^{\D}\Mt'_i$ supported in $Y_i$ is regular holonomic.
\sk
Applying the above argument to $\DD\Mt'_i$ and using the duality for the direct image of $\D$-modules by the projective morphism $\pi_i$, for instance,
$$\DD\bl(\Hc^j(\pi_i)_*^{\D}\Mt'_i\br)=\Hc^{-j}(\pi_i)_*^{\D}\1\DD\Mt'_i,$$
we can get the exact sequence
$$^{\vee}\!\!\Hc^0_{[\Zt_i|\Yt_i]}\bl(\Hc^0(\pi_i)_*^{\D}\Mt'_i\br)\buildrel{\!\!\phi_i}\over\to\Hc^0(\pi_i)_*^{\D}\Mt'_i\to\Hc^1(\pi_i)_*^{\D}\Mt'''_i.$$
Here $\Mt'''_i$ is the dual of the cokernel of the morphism
$$\DD\Mt_i\to\Hc^0_{[\Zt_i|\Yt_i]}\DD\Mt_i,$$
and is supported in $\Yt_i$. The cokernel of $\phi_i$ is the largest quotient $\D_X$-module of $\Hc^0(\pi_i)_*^{\D}\Mt'_i$ supported in $Y_i$, see Rem.\,1.2d.
The desired assertion on coherent quotients supported in $Y_i$ then follows. The assertion~(1) is thus proved.
\sk
We now show the assertion~(2) for $d_{\!\M}\eq d$ (using (1) for $d_{\!\M}\eq d\1$). Set $Z:={\rm Supp}\,\M$. Let $Y\subset Z$ be a closed subvariety such that $Z\stm Y$ is purely $d_Z$-dimensional and smooth, and the restriction of $\DRt_X(\M)[-d_Z]$ to $Z\stm Y$ is a local system. (Such a closed subvariety exists, since we have the {\it minimal\1} one satisfying the above properties by using a stratification compatible with $\DRt_X(\M)$.)
\sk
Let $Z_i$ be the irreducible components of $Z$ with dimension $d_Z$. Let $\pi'_i:\Zt_i\to Z_i$ be a desingularization such that $\Yt_i:=\pi_i^{\prime-1}(Y_i)$ is a divisor with simple normal crossings, where $Y_i\,{:=}\,Y\caps Z_i$. Let $\rho_i:\St_i\to S_i$ be a desingularization with $S_i:=f(Z_i)$. We may assume that $\ft_i:=(f|_{Z_i})\ssc\pi'_i$ factors through $\St_i$ (using a desingularization of $\Zt_i{\times}_S\St_i$).
\sk
We may assume that there is a {\it divisor\1} $\Dt_i\sst \St_i$ such that any intersection of irreducible components of $\Yt_i$ is either contained in $\ft_i^{-1}(\Dt_i)$ or smooth over $\St_i\stm \Dt_i$, replacing $S$ with an open neighborhood of a given point, and enlarging $\Dt_i$ if necessary. Indeed, the intersection of irreducible components $\Yt_{i,k}$ ($k\in J$) of $\Yt_i$ is smooth over $\St_i$ at a point of $\Yt_i$ (where $J$ is any subset of the index set of irreducible components of $\Yt_i$) if and only if we have non-vanishing at this point of the form
$$\mwdg_{k\in J}\,\ddd g_k\wedge\mwdg_{1\les j\les d_{S_i}}\ft_i^*\ddd t_j,$$
where the $g_k$ are local defining functions of $\Yt_{i,k}\sst\Zt_i$, and the $t_j$ are local coordinates of $\St_i$. This shows that the closure of the relative non-normal crossing locus is a closed subvariety. (We employ also the Bertini theorem.) Note that $d_{S_i}$ may depend on $i$.
\sk
We may assume that $\Dt_i\eq\rho_i^{-1}(D_i)$ with $D_i:=\rho(\Dt_i)$ (enlarging $\Dt_i$), and moreover there is a divisor $D\sst S$ such that $D_i\eq S_i\caps D$ (replacing $S$ with an open neighborhood of a given point, and enlarging $D_i,\Dt_i$ if necessary).
We may further assume that $f^{-1}(D)\sst Y$ (enlarging $Y_i,\Yt_i$, and taking an embedded resolution of the enlarged divisor).
\sk
From the exact sequence: $0\tos\Hc^0_{[Y]}\M\tos\M\tos\Hc^0_{[X|Y]}\M\tos\Hc^1_{[Y]}\M,$ we can deduce the short exact sequences
$$0\to\Hc^0_{[Y]}\M\to\M\to\M'\to 0,$$
$$0\to\M'\to\Hc^0_{[X|Y]}\M\to\M''\to 0,$$
with ${\rm Supp}\,\M''\sst Y$. Here $\M',\M''$ are defined by the two short exact sequences. Note that $\M''$ is regular holonomic as is proved in the proof of (1). Using the associated long exact sequence together with the inductive hypothesis and (1) for $d_{\!\M}\eq d$ as well as (2.1.1), we see that the assertion~(2) for $d_{\!\M}\eq d$ is reduced to the following:
$$\h{The $\Hc^jf^{\D}_*(\Hc^0(\pi_i)_*^{\D}\Mt_i)$ are regular holonomic ($j\ins\Z$).}
\leqno(2.1.3)$$
Set $\M_i:=\Hc^0(\pi_i)_*^{\D}\Mt_i$. We have for $j\ins\Z$
$$\aligned\Hc^j(\pi_i)_*^{\D}\Mt_i&=\Hc^j\RR\Ga_{[X|Y]}(\pi_i)_*^{\D}\Mt_i\\&=\Hc^j\RR\Ga_{[X|Y]}\Hc^0(\pi_i)_*^{\D}\Mt_i=\Hc^j_{[X|Y]}\M_i,\endaligned$$
by Prop.\,1.2 and Rem.\,1.2a (since ${\rm Supp}\,\Hc^j(\pi_i)_*^{\D}\Mt_i\sst Y$ for $j\ne 0$). Here the $\M_i$ are regular holonomic, see Rem.\,2.1a.
\sk
Using Cor.\,2.1a\,(1) and (2), Lem.\,1.2b together with the assertion~(1) for $d_{\!\M}\less d$, we then get that the $\Hc^j_{[X|Y]}\M_i$ are also regular holonomic. These imply that $\Hc^0(\pi_i)_*^{\D}\Mt_i$ in (2.1.3) can be replaced with $(\pi_i)_*^{\D}\Mt_i$ using the inductive hypothesis for (2) together with Cor.\,2.1a\,(1) and the assertion~(1) for $d_{\!\M}\less d$, since
$$\dim{\rm Supp}\,\Hc^j(\pi_i)_*^{\D}\Mt_i<d\q(j\ne0).$$
Let $\ft'_i:\Zt_i\to\St_i$ be the morphism factorizing $\ft_i:\Zt_i\to S$. By the definition of regular holonomic $\D$-modules in the beginning of {\bf 2.1}, the assertion~(2.1.3) is then further reduced to the following:
$$\h{The $\Hc^j(\ft'_i)^{\D}_*\Mt_i$ are regular holonomic,}
\leqno(2.1.4)$$
see also Rem.\,2.1a and Prop.\,1.4a--b. In the algebraic case, we may assume $S$ proper, replacing $X,S$ appropriately, for the proof of (2.1.4). The assertion~(2) for $d_{\!\M}\eq d$ then follows from Cor.\,1.5c. This finishes the proof of Prop.\,2.2.
\msn
{\bf Corollary~2.2.} (1) {\it The categories $\Msf_{\rh}(\D_X)$ are stable by the functors $\Hc^j_{[Y]}$, $\Hc^j_{[X|Y]}$ $(j\ins\Z)$ for any closed subvarieties $Y\sst X$ and by the pullbacks $\Hc^j\!f_{\D}^!$ $(j\ins\Z)$ for any morphisms $f$.}
\sk
(2) {\it Regular holonomic $\D$-modules are stable by the direct images $\Hc^j\!f^{\D}_*$ $(j\ins\Z)$ for any morphisms $f$ in the algebraic case.}
\sk
(3) {\it Regular holonomic $\D$-modules are stable by the dual functor $\DD$.}
\msn
{\it Proof.} (1) Since the assertion is local, the stability by $\Hc^j_{[Y]}$, $\Hc^j_{[X|Y]}$ follows from \h{Cor.\,2.1\1(2)} together with Lem.\,1.2b, Cor.\,2.1a and Prop.\,2.2\1(1). For the stability by $\Hc^j\!f_{\D}^!$, we may assume that $f$ is a closed immersion, since the smooth case can be shown easily by induction on $\dim{\rm Supp}\,\M$. Then the assertion follows from Rem.\,A.7a below.
\sk
(2) We can factorize $f=\fb\ssc j:X\tos Y$ so that $\fb:\Xb\tos Y$ is proper and $j:X\into\Xb$ is an affine open immersion, that is, $\Xb\stm X$ is a divisor. By Prop.\,2.2\1(2), it is enough to show the stability by $j_*$ which coincides with the sheaf-theoretic direct image in the algebraic case. Then the assertion is easily shown by induction on $\dim{\rm Supp}\,\M$ reducing to the meromorphic Deligne extension case.
\sk
(3) It is enough to show that $\DD\M$ is regular holonomic in the case $\M$ is the minimal regular holonomic extension of a local system, see Cor.\,2.1d. Set $Z:={\rm Supp}\,\M$. Let $Y\sst Z$ be the smallest closed subset such that $L:=\DRt_X(\M)[-d_Z]|_{Z\setminus Y}$ is a local system. Let $\pi:(\Zt,\Yt)\tos(Z,Y)$ be an embedded resolution inducing an isomorphism over $Z\stm Y$. There is a regular meromorphic connection $\Mt$ on $\Zt$ with poles along $\Yt$ such that $\M$ is a subquotient of $\Hc^0\pi_*^{\D}\Mt$. Then the assertion follows from the stability of regular holonomic $\D$-modules of normal crossing type under the dual functor $\DD$ (see Prop.\.1.3a) together with the commutativity of $\Hc^0\pi_*^{\D}$ and $\DD$, see Rem.\,A.4e below. This finishes the proof of Cor.\,2.2.
\msn
{\bf Remark~2.2.} For a regular holonomic $\D_X$-module $\M$, it is easy to see that
$$\DRt_X(\M)\in\Dsf^b_c(X^{\an},\C)^{[0]},$$
by reducing to the case of regular meromorphic connections by induction on $\dim{\rm Supp}\,\M$, where $\Dsf^b_c(X^{\an},\C)^{[0]}$ is the abelian full subcategory of $\Dsf^b_c(X^{\an},\C)$ constructed in \cite{BBD} (see also \cite[1.1]{lcd}). Indeed, we get two long exact sequences applying the cohomological functor ${}^{\bf p}\!\!\1\Hc^{\ssb}\DRt_X$ to the two short exact sequences written several lines before (2.1.3) where $Y$ is assumed to be a divisor. We can then verify that $\DRt_X(\M)\ins\Dsf^b_c(X^{\an},\C)^{\ges 0}$. So the commutativity of $\DRt_X$ with $\DD$ implies the desired assertion.
\bs\bs
\vbox{\centerline{\bf Appendix.~Elementary $\D$-module theory}
\bsn
In this Appendix we recall some elementary part of $\D$-module theory.}
\msn
{\bf A.1.~Left and right $\D$-modules.} There is an equivalence between the categories of left and right $\D$-modules on a smooth variety $X$. For a left $\D_X$-module $\M$, the corresponding right $\D_X$-module $\M^r$ is defined by
$$\M^r:=\om_X\sot_{\OX}\M\q\q(\om_X:=\Om_X^{d_X}),
\leqno{\rm(A.1.1)}$$
where the action of a vector field $\xi$ is given by
$$(\eta\sot m)\xi:=\eta\xi\sot m\mi\eta\sot\xi m\q\q(\eta\ins\om_X,\,m\ins\M).
\leqno{\rm(A.1.2)}$$
(Recall that the action of $\xi$ on $\om_X$ is defined by $-L_{\xi}\eta$ for $\eta\ins\om_X$ with $L_{\xi}$ the Lie derivation.) This construction can be generalized to the case of tensor product of any right $\D_X$-module and left $\D_X$-module, for instance $\D_X\sot_{\OX}\M$.
\msn
{\bf A.2.~Induced $\D$-modules and differential complexes.} For $\OX$-modules $\Lc,\Lc'$, we have the injection
$$\Hom_{\D_X}(\Lc\sot_{\OX}\D_X,\Lc'\sot_{\OX}\D_X)\into\Hom\1_{\C_X}(\Lc,\Lc'),
\leqno{\rm(A.2.1)}$$
induced by the tensor product $\sot_{\D_X}\OX$, see \cite[2.2.2]{mhp}. This can be defined also by using the isomorphism
$$\Hom_{\D_X}(\Lc\sot_{\OX}\D_X,\Lc'\sot_{\OX}\D_X)=\Hom_{\OX}(\Lc,\Lc'\sot_{\OX}\D_X),
\leqno{\rm(A.2.2)}$$
together with the first projection
$$\D_X=\OX\oplus\D_X\Theta_X\onto\OX\q(\h{induced by}\,\,\,P\mapsto P1),
\leqno{\rm(A.2.3)}$$
with $\D_X\Theta_X$ the kernel of the projection. The image of (A.2.1) will be denoted by
$$\Hom_{\Diff}(\Lc,\Lc').$$
An element of it is called a {\it differential morphism\1} from $\Lc$ to $\Lc'$.
\sk
There is a filtration $F$ on $\Hom_{\Diff}(\Lc,\Lc')$ induced by $F$ on $\D_X$ via the isomorphism (A.2.2). This is not necessarily exhaustive, for instance, if $\Lc\eq\D_X$, $\Lc'\eq\OX$ (considering the identity on $\D_X$). Set
$$\Hom_{\Diff}^f(\Lc,\Lc'):=\mcup_{p\in\N}\,F_p\1\Hom_{\Diff}(\Lc,\Lc').$$
An element of it is called a differential morphism of finite order. Put
$$\Hom_{\D_X}^f(\Lc\sot_{\OX}\D_X,\Lc'\sot_{\OX}\D_X):=\Hom_{\OX}(\Lc,\Lc')\sot_{\OX}\D_X,$$
which is identified with a subsheaf of $\Hom_{\D_X}(\Lc\sot_{\OX}\D_X,\Lc'\sot_{\OX}\D_X)$. (This is an \h{\it induced\1} $\D$-module, and is naturally isomorphic to $\Hom_{\Diff}^f(\Lc,\Lc')$, see also \cite{ind}, \cite{MS}, etc.)
\sk
We denote by $\Msf(\OX,\Diff)$ the additive category whose objects are $\OX$-modules and whose morphisms are differential morphisms in the above sense. Using the {\it injectivity\1} of (A.2.1), there is a natural functor
$$\aligned&\DR^{-1}_X:\Msf(\OX,\Diff)\to \Msf(\D_X)^r,\\&\h{with}\q\q\DR_X^{-1}(\Lc):=\Lc\sot_{\OX}\D_X.\endaligned$$
The latter module is called the {\it induced\1} $\D$-module, and $\Msf(\D_X)^r$ denotes the category of right $\D_X$-modules (similarly for $\Dsf^b(\D_X)^r$, etc.) The induced $\D$-modules are naturally defined as right $\D_X$-modules, since the de Rham functor $\DR_X$ for right $\D_X$-modules can be defined as the derived tensor product with $\OX$ over $\D_X$, which is just the inverse of $\DR_X^{-1}$ for induced $\D$-modules.
\sk
We can define the bounded complex and homotopy categories $C^b(\OX,\Diff)$, $K^b(\OX,\Diff)$ as usual. An object of these categories are called a bounded {\it differential complex.} The bounded derived category $\Dsf^b(\OX,\Diff)$ is defined by inverting $\D$-quasi-isomorphisms. Here a morphism between bounded differential complexes is called {\it $\D$-quasi-isomorphism\1} if it is a quasi-isomorphism after applying the functor $\DR^{-1}_X$ (similarly for $\D$-acyclic complexes). We can also define $D^+(\OX,\Diff)$, etc.
\sk
We denote by $\Dsf^b_{\rm coh}(\OX,\Diff)$ the full subcategory of $\Dsf^b(\OX,\Diff)$ consisting of $\Lc\ub$ with $\DR_X^{-1}(\Lc\ub)\ins \Dsf^b_{\rm coh}(\D_X)^r$.
\msn
{\bf Remark{\hs}A.2a.} There is an equivalence of categories
$$\DR^{-1}_X:\Dsf^b(\OX,\Diff)\simto \Dsf^b(\D_X)^r,
\leqno{\rm(A.2.4)}$$
with quasi-inverse given by the de Rham functor $\DR_X$, see \cite[1.5]{ind}. For a right $\D_X$-module $\M^r$, $\DR_X(\M^r)$ is defined by a complex whose component of degree $j$ is
$$\M^r\sot_{\OX}\1\mwdg^{-j}\Theta_X.
\leqno{\rm(A.2.5)}$$
For $\M^{r\1\ssb}\ins \Dsf^b(\D_X)^r$, there is a canonical quasi-isomorphism
$$\DR^{-1}_X\DR_X(\M^{r\1\ssb})\simto\M^{r\1\ssb},
\leqno{\rm(A.2.6)}$$
which implies the $\D$-quasi-isomorphism for $\Lc\ub\ins \Dsf^b(\OX,\Diff))$
$$\DR_X\DR^{-1}_X(\Lc\ub)\simto\Lc\ub,
\leqno{\rm(A.2.7)}$$
setting $\M^{r\1\ssb}:=\DR_X^{-1}(\Lc\ub)$.
\msn
{\bf Remark{\hs}A.2b.} Let $D^{+(b)}(\D_X)^r$ be the full subcategory of $D^+(\D_X)^r$ consisting of $\M^r{}\ub$ with $\Hc^j\M^r{}\ub\eq0$ ($j\gg 0$). Let $D^{+(b)}(\OX,\Diff)$ be the full subcategory of $D^+(\OX,\Diff)$ consisting of $\Lc\ub$ with $\DR^{-1}(\Lc\ub)\ins D^{+(b)}(\D_X)^r$. We have the equivalences of categories
$$\aligned \Dsf^b(\D_X)^r&\simto D^{+(b)}(\D_X)^r,\\ \Dsf^b(\OX,\Diff)&\simto D^{+(b)}(\OX,\Diff).\endaligned
\leqno{\rm(A.2.8)}$$
Indeed, the first equivalence is proved by using the truncations $\tau_{\les k}$ for $k\gg 0$.
To show the second, we can employ $\DR_X\ssc\tau_{\les k}\ssc\DR_X^{-1}$ ($k\gg 0$).
\msn
{\bf A.3.~Dual of $\D$-modules.} Let $\om_X\simto\IX\ub$ be an injective resolution as right $\D_X$-modules. (Here we do not have to assume that this is bounded using (A.2.8).) Assume
$$\Lc\ub\eq\DR_X(\M^r{}\ub)\ins \Dsf^b_{\rm coh}(\OX,\Diff)\q\h{with}\q\M^r{}\ub\ins \Dsf^b_{\rm coh}(\D_X)^r.$$
Here its important properties are as follows: The differential of the complex $\Lc\ub$ is given by differential morphisms {\it of finite order,} and there is a functorial finite filtration on $\Lc\ub$ whose graded pieces are complexes with $\OX$-{\it linear\1} differentials, that is, differential morphisms of order 0, so that $\D$-acyclicity of $\Lc\ub$ is equivalent to acyclicity of all the graded pieces.
\sk
Using (A.2.8), the dual $\DD(\Lc\ub)\ins \Dsf^b_{\rm coh}(\OX,\Diff)$ can be defined by
$$\aligned&\DD(\Lc\ub):=\Hom_{\OX}(\Lc\ub,\IX\ub)\q\q\q\h{with}\\&\DR_X^{-1}\DD(\Lc\ub)=\DD(\M^r{}\ub)\\&=\Hom_{\D_X}^f(\Lc\ub\sot_{\OX}\D_X,\I_X\ub\sot_{\OX}\D_X).\endaligned
\leqno{\rm(A.3.1)}$$
Here it is better to write $\DR_X^{-1}(\Lc\ub)$ rather than $\Lc\ub\sot_{\OX}\D_X$, since the differential of $\Lc\ub$ is {\it not\1} $\OX$-linear. The well-definedness of $\DD(\Lc\ub)$ (that is, the stability of $\D$-acyclic complexes by this functor) can be verified using the hypothesis that $\Lc\ub\eq\DR_X(\M^r{}\ub)$ in $C^b(\OX,\Diff)$.
\sk
For $\M^r{}\ub\ins \Dsf^b_{\rm coh}(\D_X)^r$, we can construct the canonical isomorphism
$$\DR^{-1}_X\DR_X(\M^r{}\ub)\simto\DD\bl(\DD(\M^r{}\ub)\br),
\leqno{\rm(A.3.2)}$$
using an involution $\iota_{\D}$ of $\I_X^j\sot_{\OX}\D_X$ exchanging the two right $\D_X$-module structures. Taking local coordinates $z_1,\dots,z_n$ with $n\eq d_X$, this involution sends $\eta\sot\mprod_i\dd_{z_i}^{\nu_i}$ (with $\eta\ins\I_X^j$) to
$$\msum_{k_1=0}^{\nu_1}\cdots\msum_{k_n=0}^{\nu_n}\,\eta\,\mprod_i\dd_{z_i}^{\nu_i-k_i}\sot\tbinom{\nu_i}{k_i}\mprod_i(-\dd_{z_i})^{k_i}.
\leqno{\rm(A.3.3)}$$
We can verify that $\iota_{\D}^2\eq{\rm id}$ using the equality $\bl(x_i\mi(x_i\mi y_i)\br)^{\nu_i}=y_i^{\nu_i}$ in the polynomial ring $\C[x_i,y_i]$, where $x_i$ and $y_i$ send $\eta\sot P$ to $\eta\dd_{z_i}\sot P$ and $\eta\sot\dd_{z_i}P$ respectively ($P\ins\C[\dd_{z_1},\dots,\dd_{z_n}]$). Here it is not necessary to use $\DR_X^{-1}\DR_X$ when we apply the second $\DD$.
\msn
{\bf A.4.~Direct image and pullback functors.} For a morphism of smooth complex varieties $f:X\tos Y$, put
$$\D_{X\to Y}:=\OX\sot_{f^{-1}\OO_Y}f^{-1}\D_Y,$$
which has a structure of a left $\D_X$ and right $f^{-1}\D_Y$-bimodule. Here the right $\D_X$-module structure is defined by using the canonical morphism $\Theta_X\tos f^*_{\!\OO}\Theta_Y$ which is the dual of the pullback $f^*_{\!\OO}\Om_Y^1\tos\Om_X^1$, where $f^*_{\!\OO}$ denotes the pullback as an $\OO$-module.
\sk
Changing the right and left $\D$-module structure, set
$$\D_{Y\leftarrow X}:=\om_X\sot_{f^{-1}\OO_Y}f^{-1}(\D_Y\sot_{\OO_Y}\om_Y^{\vee}),$$
which is a left $f^{-1}\D_Y$ and right $\D_X$-bimodule.
\sk
For $\M\ub\ins \Dsf^b(\D_X)$, $\M'{}\ub\ins \Dsf^b(\D_Y)$, their direct image and pullback are defined respectively by
$$\aligned f_*\uD\M\ub&:=\RR f_*(\D_{Y\leftarrow X}\sot^{\bf L}_{\D_X}\M\ub),\\ f^!_{\D}\M'{}\ub&:=\D_{X\to Y}\sot^{\bf L}_{f^{-1}\D_Y}f^{-1}\M'{}\ub[d_X{-}d_Y].\endaligned
\leqno{\rm(A.4.1)}$$
We employ the Godement canonical flasque resolution \cite{God} to define $\RR f_*$. For right $\D$-modules, we can define these functors by exchanging $\D_{Y\leftarrow X}$ and $\D_{Y\to X}$.
\msn
{\bf Remark{\hs}A.4a.} A left $\D$-module can be identified with an $\OO$-module endowed with an integrable connection. For a bounded complex of left $\D_Y$-modules $\M'{}\ub$, we have the canonical isomorphism
$$f^!_{\D}\M'{}\ub=\OX\sot^{\bf L}_{f^{-1}\OO_Y}f^{-1}\M'{}\ub[d_X{-}d_Y].
\leqno{\rm(A.4.2)}$$
The right-hand side has the pullback of the integrable connection corresponding to the left $\D_Y$-module structure of each component of $\M'{}\ub$ after taking a flat resolution. (Note that a flat left $\D_Y$-module $\Nc$ is flat over $\OO_Y$, since $\F\sot_{\OO_Y}\Nc\eq(\F\sot_{\OO_Y}\D_Y)\sot_{\D_Y}\Nc$ for any $\OO_Y$-module $\F$, and $\D_Y$ is flat over $\OO_Y$.)
\sk
Sometimes $f_{\D}^!$ is denoted by ${\bf L}f_{\!\OO}^*$ up to a shift of complexes. It can be shown that, in the regular holonomic case, it corresponds to the pullback functor $f^*_{\C}$ of bounded complex of $\C$-modules with constructible cohomology sheaves via the functor Sol, where the latter is the ``dual" of DR, that is, ${\rm Sol}\eq\DR\ssc\DD$.
\sk
Notice that {\it coherent\1} $\D$-modules are {\it not\1} stable by the pullback under closed immersions, for instance, if $\M'\eq\D_Y$.
\msn
{\bf Remark{\hs}A.4b.} If $f$ is a closed immersion $i:X\into Y$ and $X$ is locally defined by $z_i\eq0$ ($i\ins[1,r]$) with $z_i$ local coordinates of $Y$, then for a left $\D_X$-module $\M$, we have locally a noncanonical isomorphism
$$i_*\uD\M\cong\M[\dd_{z_1},\dots,\dd_{z_r}],
\leqno{\rm(A.4.3)}$$
where the 0-extension is omitted to simplify the notation.
\msn
{\bf Remark{\hs}A.4c.} If $f$ is a projection $X\eq Y{\times}Z\to Y$, then the direct image of a left $\D_X$-module $\M$ is defined by using the relative de Rham complex $\DR_{X/Y}(\M)$ (which is shifted by $d_Z$ to the left). This does not work if we assume only that $f$ is smooth, since there is no canonical lifting of vector fields on $Y$ to $X$. (In this case we could define only a cohomological direct image functor as in the case of Katz--Oda for Gauss-Manin connections.)
\msn
{\bf Remark{\hs}A.4d.} Assume $f$ is {\it proper.} For an induced $\D_X$-module $\M\eq\Lc\sot_{\OX}\D_X$, the direct image can be given by
$$f\uD_*\M=\RR f_*\Lc\sot_{\OO_Y}\D_Y.
\leqno{\rm(A.4.4)}$$
Indeed, proper direct image commutes with infinite direct sum, and
$$(\Lc\sot_{\OX}\D_X)\sot^{\bf L}_{\D_X}\bl(\OX\sot_{f^{-1}\OO_Y}f^{-1}(\D_Y)\br)=\Lc\sot_{f^{-1}\OO_Y}f^{-1}(\D_Y),$$
using locally a free resolution of $\Lc$. Differential morphisms are stable by the direct image $\RR f_*$ (using the Godement canonical flasque resolution), and so are differential complexes, see \cite{ind}. We have the direct image functor
$$\RR f_*:\Dsf^b(\OX,\Diff)\to \Dsf^b(\OO_Y,\Diff),
\leqno{\rm(A.4.5)}$$
together with the canonical isomorphism for $\Lc\ub\ins \Dsf^b(\OX,\Diff)$
$$\DR^{-1}_Y(\RR f_*\Lc\ub)=f\uD_*\1\DR^{-1}_X(\Lc\ub).
\leqno{\rm(A.4.6)}$$
\msn
{\bf Remark{\hs}A.4e.} Let $f:X\tos Y$ be a projective morphism of smooth complex varieties. For a coherent $\D_X$-module $\M$ having a coherent filtration $F$ locally on $Y$, there are duality isomorphisms
$$\Hc^jf_*^{\D}\1\DD(\M)\eq\DD(\Hc^{-j}f_*^{\D}\M)\q\q(j\ins\Z).
\leqno{\rm(A.4.7)}$$
This can be shown locally on $Y$ taking a locally free resolution of $\M$ (since $f$ is assumed to be projective).
\msn
{\bf A.5.~Coherent $\D$-modules.} It is known that a left $\D_X$-module $\M$ is {\it coherent\1} if and only if it is locally finitely generated over $\D_X$, and is quasi-coherent over $\OX$ (see Convention~2 at the end of the introduction).
\sk
This equivalence is rather easy to show in the algebraic case, since coherent $\OO$-modules on an affine variety ${\rm Spec}\,A$ correspond to finite $A$-modules. So we consider the analytic case. The argument is similar and simpler in the algebraic case, since we do not have to consider the closure $\Ub_x$ of an affine open neighborhood $U_x$ which is quasi-compact.
\sk
Assume the above two conditions are satisfied. There is a short exact sequence of $\D_{U_x}$-modules
$$0\to\K\to\D_{U_x}\sot_{\OO_{U_x}}\F\to\M|_{U_x}\to0,
\leqno{\rm(A.5.1)}$$
where $U_x$ is a sufficiently small open neighborhood of each point $x\ins X$ with $\F$ a free $\OO_{U_x}$-module of finite length. We may assume that the closure of $U_x$ is compact, and the above short exact sequence and the filtration $G$ are defined on the closure $\Ub_x$. This means that they are defined on some open neighborhood of $\Ub_x$. We will mean by a $\D_{\Ub_x}$-module an {\it inductive limit\1} of $\D_V$-modules with $V$ open neighborhoods of $\Ub_x$.
\sk
The filtration $F$ on $\D_{\Ub_x}$ gives the induced and quotient filtrations on $\K$, $\M|_{\Ub_x}$ respectively such that the short exact sequence is strictly compatible with $F$. For each $p\ins\Z$, we have the inclusion
$$F_p\M|_{\Ub_x}\subset G_{k(p)}\M|_{\Ub_x}\q\h{for some}\,\,\,k(p)\ins\Z,
\leqno{\rm(A.5.2)}$$
since $F_p\M|_{\Ub_x}$ is finitely generated over $\OO_{\Ub_x}$. So the $F_p\K$ are coherent $\OO_{\Ub_x}$-modules ($p\ins\Z$) by the exact sequence (A.5.1).
\sk
We have an increasing sequence of coherent $\D_{\Ub_x}$-submodules $\K_{(p)}$ of $\D_{\Ub_x}\sot_{\OO_{\Ub_x}}\F$ which are generated over $\D_{\Ub_x}$ by $F_p\K$ ($p\ins\Z$). The graded quotients $\Gr^F_{\ssb}\K_{(p)}$ form an increasing sequence of coherent $\Gr^F_{\ssb}\D_{\Ub_x}$-submodules of a free $\Gr^F_{\ssb}\D_{\Ub_x}$-module $(\Gr^F_{\ssb}\D_{\Ub_x})\sot_{\OO_{\Ub_x}}\F$. This is locally stationary considering its quotient graded modules by these graded submodules, since the supports of coherent graded $\Gr^F_{\ssb}\D_X$-modules in $T^*\!X$ are {\it conical\1} and the projection from the projectified cotangent bundle $P^*\!X$ to $X$ is proper. So the assertion follows. (Here Cartan's theorem~A is not used.)
\msn
{\bf A.6.~Coherence of direct images.} Assume $f:X\tos Y$ is {\it proper\1} and a coherent left $\D_X$-module $\M$ contains a coherent $\OX$-submodule $\F$ generating $\M$ over $\D_X$ locally on $Y$. Then it is well-known that the $\Hc^jf_*\uD\M$ ($j\ins\Z$) are coherent $\D_Y$-modules. Indeed, $\F$ defines locally on $Y$ a coherent filtration $F$ on $\M$ by $F_p\M:=F_p\D_X\,\F$ ($p\ins\Z$), and this gives a quasi-isomorphism to $\M$ from the complex whose component of degree $j$ is defined by
$$\D_X\1\sot_{\OX}\1\mwdg^{-j}\Theta_X\1\sot_{\OX}F_{p+j}\M\q\h{if}\,\,\,p\gg0,
\leqno{\rm(A.6.1)}$$
see \cite{Ka4}, etc. For the corresponding right $\D$-module $\M^r$ in (A.1.1), we have a similar complex whose component of degree $j$ is defined by
$$F_{p+j}\M^r\sot_{\OX}\1\mwdg^{-j}\Theta_X\1\sot_{\OX}\D_X.
\leqno{\rm(A.6.2)}$$
This is a typical example of an {\it induced\1} $\D$-module, and coincides with $\DR^{-1}_X\DR_X(\M^r)$ if $p\eq{+}\infty$ (that is, by omitting $F_{p+j}$). So the desired quasi-isomorphism is reduced to that
$$\Hc^j\Gr^F_p\DR_X(\M^r)\eq0\q(j\ins\Z,\,p\gg 0),
\leqno{\rm(A.6.3)}$$
locally on $X$, since $f$ is proper. This vanishing can be shown taking locally a filtered free resolution of $(\M^r,F)$.
\sk
The coherence of the direct image now follows from the quasi-isomorphism from (A.6.1) to $\M$ using the Grauert coherence theorem (see for instance \cite{Dou}) together with (A.4.4).
\msn
{\bf A.7.~Kashiwara's equivalence.} Let $i:X\into Y$ be a closed immersion of smooth complex varieties. Let $\Msf_{\rm coh}(\D_X)$ be the category of coherent $\D_X$-modules, and $\Msf_{{\rm coh},X}(\D_Y)$ be the category of coherent $\D_Y$-modules supported on $X$. The direct image functor induces the equivalence of categories, called Kashiwara's equivalence\,:
$$i_*\uD:\Msf_{\rm coh}(\D_X)\simto\Msf_{{\rm coh},X}(\D_Y).
\leqno{\rm(A.7.1)}$$
By {\bf A.5}, this is a special case of the following equivalence of categories
$$i_*\uD:\Msf(\D_X)\simto \Msf_{[X]}(\D_Y),
\leqno{\rm(A.7.2)}$$
where the right-hand side is the full subcategory of $\Msf(\D_Y)$ consisting of $\M$ satisfying the condition $\Hc^0_{[X]}\M\eq\M$, that is, any local section of $\M$ is annihilated by a sufficiently high power of the ideal of $X$, see also {\bf 1.2}.
\sk
It is easy to see that (A.7.2) is {\it fully faithful\1} taking some local coordinates, since the assertion is \h{\it local.} Using this, we can verify that the {\it essential surjectivity\1} is also {\it local.} The assertion is then reduced inductively to the case $X\eq\{z_1\eq0\}$ with $z_1,\dots,z_{d_Y}$ coordinates of $Y$ (shrinking $Y$). For $\M\ins \Msf_{[X]}(\D_Y)$, set
$$G_k\M:={\rm Ker}\,z_1^k\sst\M\,\,\,\,\h{if}\,\,\,\,k\gess 1,\,\,\,\h{and \,0\, otherwise.}$$
This filtration is exhaustive. By definition we have the injections
$$z_1:\Gr^G_{k+1}\M\into\Gr^G_k\M\q(k\gess 1).
\leqno{\rm(A.7.3)}$$
This imply that
$$z_1\dd_{z_1}{+}\1k=0\q\h{on}\q\Gr^G_k\M\q(k\gess 1),
\leqno{\rm(A.7.4)}$$
since it holds for $k\eq 1$ using $\dd_{z_1}z_1\eq z_1\dd_{z_1}{+}\1 1$, where the stability of $G_k\M$ by $z_1\dd_{z_1}$ follows from $[z_1\dd_{z_1},z_1^k]=kz_1^k$. Then the filtration $G$ has a splitting by the $\D_X$-submodules
$$\M^{(k)}:={\rm Ker}(z_1\dd_{z_1}{+}\1k)\subset\M\q(k\gess 1).$$
(Note that a polynomial relation implies the Jordan decomposition, see Rem.\,A.7c below.)
Using $[\dd_{z_1},z_1]\eq 1$, the actions of $z_1,\dd_{z_1}$ induce morphisms between the $\D_X$-submodules
$$\dd_{z_1}:\M^{(k)}\to\M^{(k+1)},\q z_1:\M^{(k+1)}\to\M^{(k)}\q(k\gess 1),
\leqno{\rm(A.7.5)}$$
with $z_1\dd_{z_1}\eq{-k}$ on $\M^{(k)}$ and $\dd_{z_1}z_1\eq {-}k$ on $\M^{(k+1)}$ ($k\gess 1$). We thus get the isomorphism
$$\M=\mopl_{k\ges 1}\,\M^{(k)}\cong\M^{(1)}[\dd_{z_1}].
\leqno{\rm(A.7.6)}$$
(Note that $\dd_{z_1}$ depends on the other coordinates $z_i$ ($i\gess 2$).)
\msn
{\bf Remark{\hs}A.7a} ({\it Relation between local cohomology and pullbacks\1}). For a $\D_Y$-module $\M$, we have the isomorphisms (see for instance \cite{Bo}):
$$i_*^{\D}\Hc^ji_{\D}^!\M=\Hc^j_{[X]}\M.
\leqno{\rm(A.7.7)}$$
Here it is convenient to use right $\D$-modules. Let $z_1,\dots,z_{d_Y}$ be local coordinates of $Y$ such that $X\eq\mcap_{i=1}^r\,\{z_i\eq0\}$ locally, where $r:=d_Y{-}d_X$. Then the pullback functor $i_{\D}^!$ is locally expressed as the Koszul complex for the actions of $z_1,\dots,z_r$, and the assertion follows from (A.7.2) by taking an injective resolution of $\M$ as a $\D$-module. 
\msn
{\bf Remark{\hs}A.7b} ({\it Relation between characteristic varieties\1}). For a coherent $\D_X$-module $\M$, we have
$${\rm Char}\,i_*^{\D}\M=\rho^{-1}({\rm Char}\M).
\leqno{\rm(A.7.8)}$$
Here $\rho:T^*Y|_X\to T^*X$ is the canonical projection, which is the dual of the inclusion of tangent bundles $TX\into TY|_X$. This follows from the definition of $i_*^{\D}$.

\msn
{\bf Remark{\hs}A.7c} ({\it Generalized Jordan decomposition\1}). Let $\F$ be a sheaf of vector spaces over a field $K$, and $T\ins{\rm End}_K(\F)$. Assume there is a polynomial $P(t)\ins K[t]$ with $P(T)\eq0$, and there is a factorization $P(t)\eq\prod_{i=1}^rP_i(t)$ with $P_i(t),P_j(t)$ mutually prime ($i\,{\ne}\,j$). Then there are mutually orthogonal projectors $\pi_i\in{\rm End}_K(\F)$ such that
$$\aligned&\pi_i^2\eq\pi_i,\q\pi_i\ssc\pi_j\eq0\,\,\,(i\,{\ne}\,j),\q\msum_{i=1}^r\,\pi_i\eq{\rm id},\\&\pi_i\eq Q_i(T)\,\,\,\h{for some}\,\,\,Q_i(t)\ins K[t],\q P_i(T)\ssc\pi_i\eq0\,\,\,(\forall\,i).\endaligned
\leqno{\rm(A.7.9)}$$
Indeed, setting $P^{\vee}_i\!(t):=\prod_{j\ne i}P_j$, we can show by induction on $r$ that there are $R_i(t)\ins K[t]$ satisfying
$$\msum_{i=1}^r\,R_i(t)P^{\vee}_i\!(t)=1,
\leqno{\rm(A.7.10)}$$
since $K[t]$ is a principal ideal domain. Then $Q_i(t)\eq R_i(t)P^{\vee}_i\!(t)$.
\sk
This implies the usual Jordan decomposition setting $T_s\,{:=}\,\sum_{i=1}^r\al_i\pi_i$ if $P_i(t)\eq(t\mi\al_i)^{m_i}$. Here any finite-dimensionality is not supposed, and $K$ is not necessarily algebraically closed, for instance, a number field. (In the case $\F$ is an $\OX$-module with $K\sst\C$ and $T\ins{\rm End}_{\OX\!}(\F)$, we have $\pi_i\ins{\rm End}_{\OX\!}(\F)$.)

\ms
{\smaller\smaller RIMS Kyoto University, Kyoto 606-8502 Japan}
\end{document}